\documentclass{amsart}

\usepackage{amsfonts}
\usepackage{amssymb}

\newcommand\R{{\mathbf R}}
\newcommand\Z{{\mathbf Z}}

\newcommand\eps{{\varepsilon}}
\newcommand\Riem{{\operatorname{Riem}}}
\newcommand\dist{{\operatorname{dist}}}
\newcommand\Hess{{\operatorname{Hess}}}
\newcommand\Ric{{\operatorname{Ric}}}
\newcommand\Vol{{\operatorname{Vol}}}
\renewcommand\min{{\operatorname{min}}}

\theoremstyle{plain}
  \newtheorem{theorem}[subsection]{Theorem}
  \newtheorem{proposition}[subsection]{Proposition}

\theoremstyle{remark}

\theoremstyle{definition}

\include{psfig}

\begin{document}

\title[Perelman's proof of Poincar\'e]{Perelman's proof of the Poincar\'e conjecture: a nonlinear PDE perspective}

\author{Terence Tao}
\address{Department of Mathematics, UCLA, Los Angeles, CA 90024}
\email{tao@@math.ucla.edu}

\begin{abstract}
We discuss some of the key ideas of Perelman's proof of Poincar\'e's conjecture via the Hamilton program of using the Ricci flow, from the perspective of the modern theory of nonlinear partial differential equations.
\end{abstract}

\maketitle

\section{Introduction}

\subsection{Perelman's theorem}

In three remarkable papers \cite{per1}, \cite{per2}, \cite{per3} in 2003, Grisha Perelman significantly advanced the theory of the Ricci flow, which was a quasilinear PDE introduced by Richard Hamilton \cite{ham1} and then extensively studied by Hamilton and others in a program to understand the topology of manifolds.  In particular, in \cite{per2} Perelman showed that in three spatial dimensions there was a well-defined \emph{Ricci flow with surgery} procedure, which is
the Ricci flow PDE punctuated at a discrete set of times by carefully chosen surgery operations which may possibly change the topology of the manifold, with the following remarkable property:

\begin{theorem}[Perelman's theorem, finite time extinction case]\label{main}\cite{per1}, \cite{per2}, \cite{per3} Let $M = (M_0,g_0)$ be a closed smooth three-dimensional compact Riemmanian manifold whose fundamental group is a free product of finite groups and infinite simple groups\footnote{This hypothesis is in fact quite natural, because the arguments also show conversely that finite time extinction for Ricci flow is only possible when the fundamental group is a free product of finite groups and infinite simple groups.  For the purposes of proving just the Poincar\'e conjecture, though, it is possible to work entirely in the category of simply connected manifolds throughout, although this only simplifies the argument at one small point (in the proof of finite time extinction).  We thank John Morgan for clarifying this point.}.  Then the Ricci flow with surgery $t \mapsto M_t = (M_t, g_t)$ is well-defined for all $t \geq 0$, and becomes extinct in finite time (thus there exists $T_*$ such that $M_t = \emptyset$ for all $t > T_*$).  Furthermore, at all times $t \geq 0$ the manifold $M_t$ is topologically equivalent to a finite union of connected sums of quotients of $S^3$ (which in particular includes the projective space $RP^3$ as a special case), and $S^2$-bundles over $S^1$.  
\end{theorem}

Applying the conclusion of this theorem at $t=0$, and using standard topological methods to detect which of the possible connected sums are simply connected, one obtains as a corollary the \emph{Poincar\'e conjecture}: every closed smooth 
simply connected three-dimensional manifold is topologically equivalent to a sphere\footnote{In fact, an additional argument, noting that surgery cannot destroy the simple-connectedness of the manifold components, demonstrates that the Ricci flow with surgery from a simply connected manifold only involves finite unions of disconnected spheres.  We thank John Morgan for clarifying this point.}.  A more difficult version of this theorem proven in \cite{per1}, \cite{per2}, \cite{per3} for general (not necessarily simply-connected) closed smooth manifolds $M = M_0$ does not have the finite extinction time result (which is false), but through an analysis of the final state (or ``scattering state'') as $t \to +\infty$, enough control on the topology of $M_t$ at each time is obtained to prove the \emph{geometrisation conjecture} of Thurston (see \cite{cao}, \cite{kleiner} for full details, including several clarifications and alternate derivations on certain difficult parts of Perelman's papers; see also the paper \cite{shioya} which fleshes out a key step in the geometrisation argument).

Perelman's result thus completes Hamilton's program and settles two major conjectures in geometry and topology, and shows the way for further progress in these fields.  But, as one can already see from the statement of Theorem \ref{main}, a large part of Perelman's work is actually conducted\footnote{It is a remarkable fact, appreciated in the last few decades (see e.g. \cite{yau}), that the flows of nonlinear PDE are an extraordinarily effective tool for understanding geometry and topology, and in particular in placing topological objects in a geometric ``normal form''.  One way to view this is that the \emph{continuous} flows of PDE, especially when augmented with a surgery procedure, are a geometrically natural successor to the more classical \emph{discrete} combinatorial algorithms that were employed to understand topology (e.g. by taking a triangulation of a manifold and then simplifying this triangulation by a sequence of combinatorial maneuvres).  While a PDE flow is in many ways ``dumber'' than a combinatorial algorithm (though this can be partially rectified by an intelligent choice of surgery procedure), if the flow is sufficiently geometrical in nature then the flows acquire a number of deep and delicate additional properties, most notably monotonicity formulae, which sometimes have no obvious analogue for discrete combinatorial algorithms.  Indeed, the non-geometric nature of most combinatorial arguments introduces a certain amount of artificial ``noise'', thus drowning out some of the more subtle monotonicity properties and making it more difficult to ensure that any given algorithm actually terminates in finite time.}
 in the arena of nonlinear PDE.  Indeed, even if one deleted the final claim in Theorem \ref{main} (thus removing the applications of this theorem to geometry and topology), the remaining arguments (and its generalisation to non-simply-connected manifolds), which are almost 
entirely nonlinear PDE arguments, would already be (in my opinion) the most technically impressive 
and significant result in the field of nonlinear PDE in recent years; the fact that this 
PDE result also gives the Poincar\'e conjecture and the more general geometrisation conjecture makes
it (again in my opinion) the best piece of \emph{mathematics} we have seen in the last ten years.  It is truly a landmark achievement for the entire discipline.

\subsection{Purpose of these notes}

Perelman's own description \cite{per1}, \cite{per2}, \cite{per3} of the arguments, which appeared in 2003, are highly succinct, with many steps which would take over a page if written out in full often compressed to a few lines; for a few years it was not entirely clear whether the argument was complete and correct, though it was evident from the very beginning that major breakthroughs had been made regardless\footnote{In particular, it was soon realised that Perelman had discovered a new monotone quantity for Ricci flow which was both \emph{critical} (scale-invariant) and \emph{coercive} (geometry-controlling).  While this by itself does not automatically prove the Poincar\'e conjecture, it made several previously impossible-looking steps in Hamilton's program now become potentially tractable.  In comparison the global regularity problem for Navier-Stokes (another major unsolved problem in nonlinear parabolic PDE) is considered extremely difficult precisely because no such critical coercive quantity for that flow is known.}.
However, in most cases these steps were of a standard type known to the experts, and could be fleshed out with reasonable work; only in a few cases were there significant enough issues with one or two steps which required more non-trivial effort to resolve\footnote{For instance, in \cite{per2} some minor corrections to \cite{per1} were noted.}, though, as noted in \cite[\S 1]{kleiner}, all such issues could be resolved using the general methods of Hamilton and Perelman.  Very recently (in 2006), three essentially simultaneous expositions\footnote{Preliminary versions of parts of \cite{kleiner} were available as early as June 2003, and influenced all three subsequent expositions.} of Perelman's work have appeared, by Kleiner-Lott \cite{kleiner}, Morgan-Tian \cite{morgan}, and Cao-Zhu \cite{cao}.  Each of these works, which emerged from large collaborative projects, some of which were in communication with each other and with Perelman, have performed an immense service to the field by giving an complete and detailed account of the entirety of Perelman's arguments\footnote{The book \cite{morgan} deals only with the Poincar\'e conjecture aspects of Perelman's work, and not on the more lengthy version of the argument establishing the geometrisation conjecture.}, and it is now certain that Perelman's original argument was indeed essentially complete and correct in every important detail.  In particular the notes here draw heavily on these three expositions as well as on the original papers of Perelman.

In these notes, I will attempt to give a high-level overview of the proof of Theorem \ref{main},
from the perspective of nonlinear evolution equations, while also attempting a sort of ``concordance'' between Perelman's papers and the three expository monographs of Kleiner-Lott, Morgan-Tian, and Cao-Zhu.  The main point of this is to try to convey some aspects of Perelman's work on Hamilton's program to those mathematicians versed in nonlinear PDE, but are not experts in either geometry or topology.  In particular, I want to convey that Perelman's work is in many ways the epitome of a nonlinear PDE argument, in which all the major milestones in the argument are familiar PDE milestones\footnote{Indeed, the argument even goes beyond standard PDE milestones, such as classification of blowup singularities, but also moves onward to new milestones, such as describing a surgery procedure to remove singularities and continuing to control the flow beyond these singularities, which are still considered too ambitious a project for most other PDE (unless one is content with only possessing very weak solutions for which one has little control).}, but where the execution requires extremely delicate analysis and several major technical breakthroughs, many of which arise by carefully exploiting the rich geometric structure of Ricci flow. 

In order to keep things as simple as possible I will only focus on Theorem \ref{main}, which is sufficient for Poincare's conjecture, and not discuss the more complicated version of this theorem which also yields the more general geometrisation conjecture.  This simpler theorem already captures a large fraction of the breakthroughs that Perelman made, and should already be sufficient for conveying the points discussed in the previous paragraph.

Because the focus of these notes will be on the nonlinear PDE aspects of the work, the other aspects from geometry and topology will be downplayed.  In particular, various geometric structural equations arising from Ricci flow shall be presented simply as miraculous identities, with no derivation or intuitive geometric explanation given; also, the use of standard topological tools to convert one type of topological control on a manifold (such as control of homotopy groups) to another will be sketched only briefly.  Thus these notes should not be in any way construed as a balanced presentation of Perelman's work or more generally of Hamilton's program, as it shall be by conscious design heavily slanted towards the PDE component of the work.  This is in part for reasons of space and time, but also because I am not sufficiently qualified in geometry or topology to expertly comment on the other components of these papers. 

We will not attempt here to give anything resembling a complete or detailed presentation of Perelman's work, or even of just the nonlinear PDE aspects of this work; for such presentations, we refer the reader to \cite{cao}, \cite{kleiner}, \cite{morgan}.  In particular we will often choose to sacrifice rigour for sake of exposition, and in fact will commit a number of small lies\footnote{We will attempt to apologise for these lies using footnotes whenever feasible, and cross-reference to the more rigorous statements in other papers when appropriate.} in order to suppress certain technical issues which might otherwise obscure the underlying ideas.  In a similar spirit, we will often not discuss the strongest version of any given component of the proof, but instead focus on a much simpler model proposition which is not strong enough for the full proof of Poincar\'e, but does indicate the main ideas that are used in the proof of the true proposition.
We will also not attempt to systematically survey the literature and other context which Perelman's work and Hamilton's program builds upon and is placed in; for this, we refer the reader to \cite{milnor}, \cite{morgan-2}, \cite{anderson}, \cite{yau-poincare}.  Finally, it should be noted that the arguments here are arranged in an ahistorical order, in order to expedite the natural flow of the proof. In particular, recent arguments (such as those of Perelman) will sometimes be juxtaposed with earlier arguments (such as those of Hamilton), reflecting the fact that different parts of Hamilton's program to prove the Poincar\'e conjecture via Ricci flow were understood at different times.  

As has been made clear in several other places (see e.g. the introduction of \cite{cao}), the primary contributions here are those of Hamilton and Perelman.  Of course, as with all other modern mathematical achievements, many other mathematicians in the field have also made important contributions, either by laying foundational work, establishing tools, developing insights or suggestions, or by clarifying, refining, generalising, or giving alternate proofs of key results.  In particular, as has been mentioned earlier, the efforts in \cite{cao}, \cite{kleiner}, \cite{morgan} to fully explicate all the details of Perelman's arguments have been a significant, and extremely non-trivial, service to the field. Assigning credit to each particular idea or insight is of course difficult to do with absolute accuracy.  While we will try our best to attribute things here, we will refer the reader to the careful and detailed expositions in \cite{cao}, \cite{kleiner}, \cite{morgan} for a more thorough treatment of each of the steps and their provenance.

\subsection{Overview of argument}

In the jargon of nonlinear PDE, the Ricci flow is a quasilinear parabolic system, and Theorem \ref{main} is (in part) a finite time blowup result\footnote{In particular this suggests that the nonlinear component of this system is ``focusing'' in nature, acting in the direction of blowup rather than in the direction of dissipation.} from essentially arbitrary data for this system.  But it is also much more than merely a finite time blowup result; it gives quite precise control of the behaviour near each blowup point, and consequently also allows one to continue the evolution \emph{beyond} each blowup point (which is a level of understanding which is currently not enjoyed by many commonly studied nonlinear PDE, even those that are substantially simpler\footnote{This is perhaps an unfair comparison.  While the Ricci flow is much more nonlinear than some other model nonlinear equations, such as the semilinear heat, wave, and Schr\"odinger equations, it is parabolic rather than hyperbolic or dispersive, and perhaps more importantly is an extremely \emph{geometric} equation, thus enjoying significantly more structure than most other model nonlinear PDE.  The combination of the parabolic and the geometric natures of the flow in particular leads to a very large number of monotonicity formulae, which are essential to the analysis and which would not be expected to be present for many other model nonlinear PDE.  We thank Igor Rodnianski for emphasising this point.}).  The need to continue the evolution even after the appearance of singularities is essential for the topological applications of the theorem, because the topological structure of the evolved manifold $M_t$ is only obvious at late times (in particular, after the extinction time $T_*$), and then to reconstruct the topology of the initial manifold $M_0$ one needs to understand the nature of the blowup at all blowup times between $t=0$ and the extinction time $t=T_*$, and not just the first such blowup time\footnote{In the groundbreaking earlier work of Hamilton \cite{ham1}, it was shown that if the manifold $M_0$ initially has positive Ricci curvature everywhere, then there is extinction in finite time with no intervening blowups, and that the metric becomes asymptotically round.}.

The proof of Theorem \ref{main} follows standard PDE paradigms, even if the execution of these paradigms is extraordinarily non-trivial.  Firstly, a local-in-time existence result for the Ricci flow
is needed.  This was first established by Hamilton \cite{ham1}, with a simpler proof given by DeTurck 
\cite{turck}, and is sketched out in Section \ref{local-sec}.  With the local existence result in hand, Theorem \ref{main} then reduces to a structural result on the blowup times, and on the nature of the blowup singularities at each of these times.

To achieve the latter task, one applies a standard rescaling limit argument around each blowup singularity to understand the asymptotic profile of the blowup.  These asymptotic profiles were established in special cases by Hamilton \cite{hamil-sing} before Perelman's work, and conjectured in general.  One major obstacle 
to establishing this conjecture was the need for \emph{critical} (scale-invariant) and \emph{coercive} (geometry-controlling) quantities  as one approached a singularity of Ricci flow.
For instance, a paper of Hamilton \cite{hamil-sing} largely classifies the asymptotic profiles of singularities under the scale-invariant assumption that the renormalised curvatures and injectivity radius stayed bounded.  However, it was not known at that time how to attain such scale-invariant assumptions.

Previously to Perelman's work, the type of controlled quantities which were known were either \emph{supercritical} in nature, which meant that they became dramatically weaker for the purposes of controlling the solution as one approaches the singularity, or \emph{non-coercive}, thus not controlling all aspects of the geometry (except under special assumptions such as positive curvature).  This raises the possibility that the rescaled solutions around that singularity might not be sufficiently bounded to extract a limit.  The first major breakthrough in Perelman's work is to introduce a number of new \emph{critical} monotone quantities, which are coercive enough to provide adequate boundedness on the rescaled solutions (in particular, uniform bounds on the Riemann curvature) to provide a limit (by Hamilton's compactness theory \cite{hamil-compact} for Ricci flow, based in turn on the Cheeger-Gromov theory \cite{cheeger} of limits of Riemannian manifolds).  This very important \emph{no-collapsing theorem} of Perelman\footnote{Actually, Perelman establishes a variety of no-collapsing theorems with variations in the hypotheses and conclusions, and for Ricci flow with and without surgery, as do the authors following him: see \cite[\S 4,7,8]{per1}, \cite[\S 7]{per2}, \cite[\S 8]{morgan}, \cite[\S 12, 25, 27, 78]{kleiner}, \cite[\S 3.3, 3.4, 7.6]{cao}.} controls the renormalised geometry as one approaches blowup, and allows one to classify the asymptotic profiles of blowup singularities; a refinement of these arguments then allows one to understand the geometry of high-curvature regions quite accurately also.  The latter is then crucial for defining a surgery procedure properly.

The next major task is to establish finite extinction time\footnote{This is the ``shortcut'' approach to proving the Poincar\'e conjecture.  There is a more difficult approach that gives not only this conjecture, but also the geometrization conjecture, which allows for the Ricci flow-with-surgery to evolve for infinite time, but we will not discuss this more complex approach here.} for simply connected manifolds.  This was achieved simultaneously by Perelman \cite{per3} and by Colding-Minicozzi \cite{cold}, and also explained in detail in \cite[\S 18]{morgan}.  While non-trivial, these results are somewhat simpler in nature than those in \cite{per1}, \cite{per2}, relying on a different monotone quantity from that considered earlier\footnote{In particular, this new monotone quantity is no longer critical, but criticality is not as relevant for blowup time bounds as it is for blowup singularity analysis, as rescaling is not involved.}  which was known to be both decreasing and positive for Ricci flow, which can then be used to deduce finite extinction time (after carefully attending to the behaviour of these quantities during surgery).  

Finally, there is the endgame in which one returns back from the final extinction time to the original
state of the Riemannian manifold, and verifies that the topology is always a finite union of spheres.  This turns out to be relatively easy once the other two major tasks are accomplished.

\subsection{Acknowledgements and history}

The author thanks Igor Rodnianski for helpful discussions and John Lott, Gang Tian and especially John Morgan for many clarifications and insights. The material here is of course heavily based on several of the papers in the references, but the author was particularly reliant on all three of the excellently detailed expositions of Perelman's work in \cite{kleiner}, \cite{morgan}, \cite{cao}, in addition to Perelman's original papers \cite{per1}, \cite{per2}, \cite{per3}.  Indeed, while these papers clearly have strong overlap (as should be obvious from the discussion below), they each have slightly different perspectives, strengths, and weaknesses, and by reading all these papers in parallel (in particular, switching from one paper to another whenever I was stuck on a particular point) I was able to obtain a far richer view of the argument than I could have obtained just from reading one of them.

An initial version of this manuscript was posted on my web site on October 9, 2006.  After receiving several comments and corrections, a revised version was then placed on the ArXiV on October 29, 2006.
At this present time, I am not planning to submit this manuscript for publication.

\section{Notation}

I will be using notation from Riemannian geometry, which I will define very rapidly here without providing any motivation as to where these geometric definitions and identities arise from.  For reasons of personal preference I tend to express tensors using coordinates, for instance the metric $g$ would be expressed as $g_{\alpha \beta}$, with the usual raising and lowering and summation conventions.  In many other treatments it is more customary to use coordinate-free notation.

Throughout these notes, all Riemannian manifolds $M = (M,g)$ will be smooth compact and three-dimensional; while a substantial portion of Perelman's work does extend to general dimension, we will not discuss those generalisations in order to slightly simplify the discussion.

Given a Riemannian metric $g = g_{\alpha \beta}$, one has the \emph{inverse metric} $g^{-1} = g^{\alpha \beta}$, and the \emph{volume form} $dg = \sqrt{\det(g_{\alpha \beta})} dx$.  The Levi-Civita connection $\nabla$ is then defined in coordinates via the \emph{Christoffel symbol}
$$ \Gamma^\alpha_{\beta \gamma} := \frac{1}{2} g^{\alpha \delta} (\partial_\beta g_{\gamma \delta} + \partial_\gamma g_{\beta \delta} - \partial_\delta g_{\beta \gamma}).$$
This in turn leads to the \emph{Riemann curvature tensor}
$$ \Riem_{\alpha \beta \gamma}^\delta := \partial_\alpha \Gamma_{\beta \gamma}^\delta - \partial_\beta \Gamma_{\alpha \gamma}^\delta + \Gamma_{\beta \gamma}^\sigma \Gamma_{\alpha \sigma}^\delta -
\Gamma_{\alpha \gamma}^\sigma \Gamma_{\beta \gamma}^\delta$$
which then contracts to form the \emph{Ricci curvature tensor}
$$ \Ric_{\alpha \gamma} := \Riem_{\alpha \beta \gamma}^\beta$$
which in turn contracts to form the \emph{scalar curvature}
$$ R := g^{\alpha \gamma} \Ric_{\alpha \gamma}.$$
From the point of view of dimensional analysis\footnote{Very informally, one can think of this dimensional analysis as simply counting the number of times the symbols $g$ and $\partial$ appear in every expression.  Ricci flow will give $\partial_t$ the same dimension as $O( g^{-2} \partial^2 )$.  See Table \ref{dimtable}.}, it is helpful to keep in mind the schematic form of these tensors (thus omitting indices, signs, and constants):
\begin{align*}
\Gamma &= O( g^{-1} \partial g )  \\
\Riem, \Ric &= O( g^{-1} \partial^2 g ) + O( g^{-2} (\partial g) (\partial g) )  \\
R &= O( g^{-2} \partial^2 g ) + O( g^{-3} (\partial g) (\partial g) ).
\end{align*}
In three spatial dimensions, one also has $dg = O(g^{3/2}) dx$.  Of course, these schematic forms destroys all the geometric structure, which is absolutely essential for such tasks as establishing monotonicity formulae, but for the purposes of analytical aspects of the Ricci flow, such as the local existence theory, these schematic forms are often adequate.

Given two orthonormal tangent vectors $X = X^\alpha$, $Y = Y^\beta$ at a point $x$, the \emph{sectional curvature} $K(X,Y)$ of the tangent plane spanned by $X$ and $Y$ is given by the formula
$$ K(X,Y) = \Riem_{\alpha \beta \gamma}^\delta X^\alpha Y^\beta X^\gamma Y_\delta.$$
One can show that this curvature depends only on the plane spanned by $X$ and $Y$, indeed it can be interpreted as the Gauss curvature of this plane at $x$.  Schematically, $K$ has the same units as the Ricci scalar: $K = O( g^{-2} \partial^2 g ) + O( g^{-3} (\partial g) (\partial g) )$.  Indeed one can think of $R$ (up to an absolute constant) as the expected sectional curvature of a randomly selected $2$-plane, whereas the Ricci bilinear form $\Ric_{\alpha \beta} X^\alpha X^\beta$ for a unit tangent vector $X$ is the expected sectional curvature of a randomly selected $2$-plane passing through $X$.

The \emph{Laplace-Beltrami} operator $\Delta_g$ is defined on scalar functions $f$ as
$$ \Delta_g f = g^{\alpha \beta} \partial_\alpha \partial_\beta f - g^{\alpha \beta} \Gamma_{\alpha \beta}^\gamma \partial_\gamma f $$
or schematically as
$$ \Delta_g f = O( g^{-1} \partial^2 f ) + O( g^{-2} (\partial g) (\partial f) ).$$
It is worth noting that $\Delta_g$ has the same scaling as the Ricci scalar $R$.  As such, the Ricci scalar often naturally arises as a lower order correction term to the Laplacian.  This is for instance evident in the formula for the \emph{conformal Laplacian} $\Delta_g + \frac{1}{8} R$ from conformal geometry, although this Laplacian does not seem to directly play a role in Perelman's work\footnote{On the other hand, to verify various monotonicity formulae here one does need to apply the Laplace-Beltrami operator to forms as well as to scalar fields, and the Ricci tensor again arises as a lower order term in those cases.  However as we are not pursing the structural geometric identities in detail here we shall not go into these issues.}.

The Ricci flow, introduced for the study of topology of manifolds by Hamilton \cite{ham1}, 
is a flow which (assuming a suitable local existence theory) 
produces a curve $t \mapsto (M_t, g_t)$ of three-dimensional Riemannian manifolds for all time $t$ in some time interval. It is usually described by keeping the manifold $M_t$ fixed (thus $M_t = M$, though perhaps $M_t = \{t\} \times M$ would be more appropriate) and varying the metric $g = g_t$ via the partial differential equation
$$ \partial_t g_{\alpha \beta} = - 2 \Ric_{\alpha \beta}.$$
The factor $2$ is a traditional and convenient normalisation which is not of essential importance, however the minus sign $-$ is crucial in order for Ricci flow to be well-posed in the forward time direction.  Intuitively, Ricci flow expands (and flattens) the negative curvature portion of the manifold while contracting the positive curvature portion.

In order to continue Ricci flow past a singularity, it is often convenient to allow for a more general coordinatisation of the flow, in which the three-dimensional ``spatial'' manifolds $M_t$ are viewed as level sets of a time function $t$ on a four-dimensional ``spacetime'' manifold (possibly containing some singularities or boundaries caused by surgery), and the time derivative $\partial_t$ is replaced by a Lie derivative ${\mathcal L}_T$ along some time vector field $T^\alpha$ (with the normalisation ${\mathcal L}_T t = 1$).  See \cite[\S 4]{per2}, \cite[\S 3.8, \S 14]{morgan}, \cite[\S 67]{kleiner}, \cite[\S 7.3]{cao} for a formalisation of these ``generalised Ricci flows'', as well as the related (but subtly different) ``Ricci flows with surgery''.  To simplify the exposition we shall not discuss these flows (though they are of course essential in the full proof of the Poincar\'e conjecture) and instead speak (somewhat imprecisely) solely about ordinary Ricci flow (in which the manifold $M$, and in particular its topology, remains invariant) together with isolated surgery operations.  Nevertheless it is still useful to keep in mind the idea that at two different times $t_1, t_2$, the manifolds $M_{t_1}$ and $M_{t_2}$ are really distinct manifolds, even if they can in principle be identified with a single manifold $M$ (at least if there is no surgery operation intervening between $t_1$ and $t_2$).  To emphasise this, we shall sometimes refer to a point $x$ on $M_t$ as $(t,x)$ rather than just $x$, thus we view points as elements of spacetime rather than merely of space.

In order to avoid discussing a lot of technicalities regarding how various small parameters depend on other small parameters (which are actually quite important to the detailed analysis, but which are too distracting for a high-level overview) we shall adopt the following rather imprecise notation.  We use $X = O(Y)$, $X \lesssim Y$, or $Y \gtrsim X$ to denote the assertion that $X$ is dominated by a bounded multiple of $Y$; this bound may in fact depend on a number of minor parameters but should at least be independent of things like the natural spatial scale or the time until blowup.  Similarly we use $X \sim Y$ to denote $X \lesssim Y \lesssim X$.

\section{Dimensional analysis}

In the large data global theory of nonlinear evolution equations, the role of scale invariance, and the subsequent distinction of conserved or monotone quantities as \emph{subcritical}, \emph{critical}, and \emph{supercritical}, plays a fundamental role.  Indeed, it is not too much of an exaggeration to say that with current technology, a large data global existence result is only possible either in the presence of a critical or subcritical controlled quantity (by which we mean a conserved or monotone quantity).  Furthermore, global \emph{asymptotics} and \emph{blowup profiles} 
are only possible in the presence\footnote{Furthermore, these quantities need to be \emph{coercive}, in that bounds on such quantities severely restrict the possible range of solution behaviours.  For example, a quantity which could be bounded due to a cancellation between a large positive and large negative term is unlikely to be coercive, unless control from other sources can somehow eliminate this cancellation scenario.}  of a critical controlled quantity, or the combination of a subcritical and a supercritical controlled quantity.

In the context of Ricci flow, the scale invariance is given in coordinates by
\begin{equation}\label{gab}
 g_{\alpha \beta}(t,x) \mapsto \lambda^2 g_{\alpha \beta}( \frac{t}{\lambda^2}, x )
 \end{equation}
where $\lambda > 0$ is arbitrary.  Comparing this with the (local) scaling diffeomorphism
\begin{equation}\label{local-scale}
 g_{\alpha \beta}(t,x) \mapsto \frac{1}{\lambda^2} g_{\alpha \beta}( t, \frac{x}{\lambda} )
 \end{equation}
(which is just a change of coordinates, and which therefore will preserve every \emph{local}\footnote{The key word here is ``local''.  A global quantity such as the total volume $\int_M\ dg$ is geometric but not invariant under this scaling diffeomorphism unless one is willing to rescale the domain of the coordinate chart(s) of $M$ to be arbitrarily large or small.  On the other hand, a localised volume such as $\int_M e^{-\dist_{M_t}(x,y)^2 / 4(t_0-t) }\ dg$ would be asymptotically invariant under \eqref{local-scale} in the limit where $t_0-t$ is small, as this integral is effectively localised to a small, scale-invariant, region of space.} scalar geometric quantity), one can also write the scale invariance locally in the form
\begin{equation}\label{g-scale}
 g_{\alpha \beta}(t,x) \mapsto g_{\alpha \beta}(\frac{t}{\lambda^2}, \frac{x}{\lambda})
\end{equation}
which is the familiar parabolic scaling.

A \emph{critical} geometric quantity is one which is invariant under the scale invariance (being geometric, it does not matter which of the above two definitions of scaling one chooses).  Such quantities will be essential for understanding behaviour of the flow near singularities; subcritical quantities would be preferable, but unfortunately turn out to be insufficiently \emph{coercive} to control the geometry by themselves, while supercritical quantities (such as the total volume of the manifold) become increasingly useless as one approaches the singularity.  In parabolic theory, one can sometimes compensate for the lack of critical controlled quantities by relying instead on the maximum principle; this was for instance achieved by Hamilton \cite{hamil-sing} when establishing global existence (and absence of singularities) for two-dimensional Ricci flow assuming that the Ricci curvature is initially negative definite.  However it appears that the maximum principle by itself is not strong enough to control the evolution of Ricci flow in the case of indefinite curvature.

One quick informal way to determine the criticality of an expression is to perform a dimensional analysis.  It is convenient to do this by using two units, a unit $L$ for the spatial coordinate scale, and a unit $G$ for the metric $g_{\alpha \beta}$.  This allows us to compute the dimension of every quantity in these notes in terms of $G$ and $L$; see Table \ref{dimtable}.

\begin{table}[ht]
\begin{tabular}{|l|l|l|l|l|l|l|}
\hline 
        &  $L^{-2}$          &  $L^{-1}$         &  $L^0$    &   $L^1$      &    $L^2$  & $L^3$\\
\hline
$G^{-1}$&
$\left.\begin{array}{l}
\partial_t, \Delta_g,R,R_\min,\\
|\Riem|_g,|\Ric|_g, K
\end{array}\right.$
                 &$\dot \gamma$      & $g^{-1}$     &              &           &\\
\hline
$G^{-1/2}$&                  &$|\dot \gamma|_g$  &              &              &           &\\
\hline
$G^0$     &$\Riem, \Ric, \pi$&$\partial,\Gamma,X$&$\left.\begin{array}{l}
\tilde V,l,\\
\kappa,f,\delta
\end{array}\right.$&$x,x_0,\gamma$&           &\\
\hline
$G^{1/2}$ &                  &                   &              &
$\left.\begin{array}{l}\dist, {\mathcal L}, r,\\
h, \rho
\end{array}\right.$       &           & \\
\hline
$G^1$     &$\partial^2 g$    &$\partial g$       & $g$          &              &$t,t_0,\tau,\int d\tau$&\\
\hline
$G^{3/2}$ &                  &                   &           &              &           &$\Vol$,$\int dg$\\
\hline
\end{tabular}
\caption{The dimension of various quantities encountered in this paper, as powers of $G$ and $L$.  Note that scalar geometric quantities all lie on the diagonal where the power of $L$ is twice that of $G$; this is a manifestation of the dilation symmetry \eqref{local-scale}.  Tensor quantities typically do not lie on this diagonal, although their magnitudes do.  $\pi$ here refers to the deformation tensor of a vector field $X$.}
\label{dimtable}
\end{table}

It is the power of $G$ which determines the criticality, as can be seen from \eqref{gab}.  A net negative power of $G$ implies that the quantity is sub-critical (control of this quantity becomes increasingly powerful at fine scales and close to blowup times, or equivalently when the metric becomes large), whereas a positive power of $G$ implies a super-critical quantity (and hence of little use near a blowup, unless one has a favourable sign).  A quantity with no powers of $G$ is critical and thus a good candidate for establishing uniform control on rescaled solutions near a singularity (provided of course that it is sufficiently ``coercive'', in that control on the quantity implies significant control on the geometry).

The power of $L$ is linked to that of $G$ via the local scaling \eqref{local-scale}; as long as the quantity is scalar and somehow ``local'' in nature, as well as being geometric, the power of $L$ must be twice the power of $G$.  There is also a relationship for tensor quantities depending on the number of subscripts and superscripts in the tensor but it is slightly more tricky to state; it will be left to the reader (perhaps using the above table as guidance).

The reader is also invited to verify the dimensional correctness of all the equations in these notes.

\section{Local existence}\label{local-sec}

Given that we are trying to establish a large data final state result for a nonlinear PDE, it is of course
natural to begin with the local existence theory.  Here we run into an initial obstacle, which is that the Ricci flow
$$ \partial_t g_{\alpha \beta} = - 2 \Ric_{\alpha \beta}$$
does not obviously have the structure of a quasilinear parabolic equation.  Indeed, if we expand the Ricci tensor in
coordinates we obtain
\begin{equation}\label{ric-decomp}
\begin{split}
 -2\Ric_{\alpha \beta} &= g^{\gamma \delta}( \partial_{\gamma \delta} g_{\alpha \beta} + \partial_{\alpha \beta} g_{\gamma \delta} - \partial_{\beta \gamma} g_{\alpha \delta} - \partial_{\alpha \gamma} g_{\beta \delta}) + O( g^{-2} (\partial g) (\partial g) )\\
 &= O( g^{-1} \partial^2 g ) + O( g^{-2} (\partial g) (\partial g) ).
\end{split}
\end{equation}
Thus the Ricci flow equation takes the schematic form
$$ \partial_t g = O( g^{-1} \partial^2 g ) + O( g^{-2} (\partial g) (\partial g) )$$
where the lead term $O( g^{-1} \partial^2 g )$ is not elliptic.  This is not a manifestly parabolic equation.

Despite this, it is still possible to establish local existence for the Ricci flow for suitably regular initial data.  There are two approaches.  The first approach, due to Hamilton \cite{ham1}, uses the Nash-Moser iteration method to compensate for the lack of smoothing present in the equation.  The second approach, due to de Turck \cite{turck}, is a ``gauge fixing'' approach, taking advantage of the geometric nature of the Ricci flow, which in practical terms creates the ``gauge invariance'' formed by diffeomorphisms of the manifold (i.e. changes of coordinates).  We shall briefly discuss the latter method here.

Recall that every smooth vector field $X^\alpha = X^\alpha(x)$ on a manifold $(M,g)$ generates an infinitesimal diffeomorphism $x^\alpha \mapsto x^\alpha + \eps X^\alpha + O(\eps^2)$.  Pulling back the metric $g_{\alpha \beta}$ by this diffeomorphism creates a infinitesimally deformed metric $g_{\alpha \beta} \mapsto g_{\alpha \beta} + \eps \pi_{\alpha \beta} + O(\eps^2)$, where
$$ \pi_{\alpha \beta} := {\mathcal L}_X g_{\alpha \beta} = \nabla_\alpha X_\beta + \nabla_\beta X_\alpha$$
is the \emph{deformation tensor} of the vector field $X$.  Thus, for instance, if the deformation tensor vanishes, then $X$ generates an infinitesimal isometry, and is known as a Killing field.

The Ricci flow commutes with diffeomorphisms.  Because of this, the Ricci flow equation is gauge equivalent to the equation
$$ \partial_t g_{\alpha \beta} = - 2 \Ric_{\alpha \beta} + \pi_{\alpha \beta}$$
where $\pi_{\alpha \beta}$ is the deformation tensor of an \emph{arbitrary} (time-dependent) vector field $X$.  Indeed,
the solution to this gauge transformed Ricci flow is (formally) equal to the solution to the original Ricci flow with the same initial data, composed with the diffeomorphism obtained by integrating the vector field $X$; we omit the details.  In practical terms (and ignoring some technical issues in justifying uniqueness, preserving regularity, etc.) this 
means that we have the
freedom to add any expression of the form $\pi_{\alpha \beta} = \nabla_\alpha X_\beta + \nabla_\beta X_\alpha$ to the right-hand side of
the Ricci flow equation.  Note from the form of the Christoffel symbols $\Gamma = O( g^{-1} \partial g )$ that
$$ \pi_{\alpha \beta} = \partial_\alpha X_\beta + \partial_\beta X_\alpha + O( g^{-1} (\partial g) X ).$$
In particular if one takes\footnote{This is not a completely arbitrary choice of infinitesimal diffeomorphism; it turns out that the one-parameter family of diffeomorphisms generated by this recipe is also the harmonic map heat flow deformation of the identity map from $M$ to itself.  See \cite{lutian} for further discussion.}
$$ X_\alpha := g^{\gamma \delta} ( \partial_\gamma g_{\alpha \delta} - \frac{1}{2} \partial_\alpha g_{\gamma \delta} ) = O( g^{-1} \partial g )$$
then we see that
$$
\pi_{\alpha \beta} = g^{\gamma \delta}( - \partial_{\alpha \beta} g_{\gamma \delta} + \partial_{\beta \gamma} g_{\alpha \delta} + \partial_{\alpha \gamma} g_{\beta \delta}) + O( g^{-2} (\partial g) (\partial g) ).
$$
Thus this gives a gauge-transformed Ricci equation of the form
$$ \partial_t g_{\alpha \beta} = g^{\gamma \delta} \partial_{\gamma \delta} g_{\alpha \beta} + O( g^{-2} (\partial g) (\partial g) )$$
or in other words
$$ \partial_t g = \Delta_g g + O( g^{-2} (\partial g) (\partial g) ).$$
This is now more obviously a parabolic equation, being a nonlinear perturbation of the heat equation.  One can now use standard parabolic methods, such as iteration methods, to construct\footnote{One cautionary note regarding the use of gauge transforms is that uniqueness of the gauge-transformed equation does not always imply uniqueness of the original equation, because of the possible presence of exotic solutions for which the gauge transform is unavailable due to lack of regularity.  Nevertheless there is a satisfactory uniqueness theory for Ricci flow, see e.g. \cite{chenzhu}.  However one can still make the proof of the Poincar\'e conjecture work even if uniqueness was not known; it makes a few minor facts harder to prove (such as the fact that spherically symmetric data leads to spherically symmetric solutions) but is otherwise not important to the argument.} local solutions.

It is worth remarking that while de Turck's gauge transformed Ricci flow equation is very convenient for establishing local existence, it appears that for the global theory it is more useful to retain the original Ricci flow equation.  One reason for this is that the metric $g$, while of course essential in the local theory, is not really the main actor in the global theory; it is other quantities derived from the metric, such as the Ricci curvature tensor and the reduced length and volume, which play a much more dominant role.  The Ricci flow equation in its original form seems well suited for studying these latter quantities.

On the other hand, there is a specific gauge transformed version of the Ricci flow which is worth noting.  Suppose that the vector field $X$ is the negative gradient of a potential function $f$: $X^\alpha = -\partial^\alpha f$.  Then the Ricci flow equation transforms to
$$ \partial_t g_{\alpha \beta} = - 2 \Ric_{\alpha \beta} - 2 \Hess(f)_{\alpha \beta}$$
where $\Hess(f)_{\alpha \beta} = \nabla_\alpha \nabla_\beta f$ is the Hessian of $f$.  Now suppose that at time $t=-1$, the manifold obeys the \emph{gradient shrinking soliton equation}
\begin{equation}\label{gss}
\Ric_{\alpha \beta} + \Hess(f)_{\alpha \beta} - \frac{1}{2} g_{\alpha \beta} = 0
\end{equation}
then the Ricci flow extends backwards and forwards in time for all time $t \in (-\infty,0)$, with the solution at time $t$ being diffeomorphic to the solution at time $-1$, with metric scaled by $-t$.  This type of solution is known as a \emph{gradient shrinking soliton}, and plays a key role in the analysis of blowup singularities, similar to the role played by solitons in other equations (particularly equations with a critical conserved or monotone quantity).  There is also the \emph{gradient steady solitons}, which solve the equation
$$
\Ric_{\alpha \beta} + \Hess(f)_{\alpha \beta} = 0
$$
and which exist for all time and are constant in time up to diffeomorphism.

More refined analysis can give some explicit blowup criteria necessary for singularities to form.  For us, the most important one is that if the Ricci flow cannot be continued past a certain time, then the Riemann curvature must be blowing up to infinity at that time; see \cite{shi-thesis} or \cite[Theorem 6.3]{chow}.  Thus singularities of the Ricci flow will be preceded by the appearance of regions of arbitrarily large curvature.

The above local existence and uniqueness theory is not only valid for compact manifolds, but extends to non-compact manifolds of bounded curvature; see \cite{chenzhu}.  This turns out to be important later on when we construct
``standard solutions'' which are needed to replace singularities in the surgery procedure.

\section{Monotonicity formulae I: Lower bounds on curvature}

Perhaps the primary tools for obtaining long-time control of parabolic equations are monotonicity
formulae\footnote{This is in contrast with Hamiltonian equations such as wave or Schr\"odinger equations, in which the primary tool is conservation laws.  Indeed the only important conserved quantity of Ricci flow is a discrete quantity, namely the topology of the underlying manifold (though of course this conservation can be broken during surgery).  Of course, without this conserved quantity the whole program of trying to extract topological information from Ricci flow would be doomed from the start!}.  

These formulae arise from a number of sources.
Prior to Perelman's work, it was already known that \emph{maximum principles} (including the classical scalar maximum principle, but also the more advanced tensor maximum principle of Hamilton) yielded some very useful monotonicity formulas which controlled the curvature from below, but unfortunately not from above - so the monotone quantities were not fully \emph{coercive}.

For instance, the minimal scalar curvature $R_\min(t) := \inf_{x \in M_t} R(t,x)$ was a non-decreasing function of $t$; this follows easily from the observation that this scalar curvature obeys a nonlinear heat equation\footnote{For proofs of the identities in this section, see e.g. \cite{chow}, \cite{cao}.}
\begin{equation}\label{r-eq} \partial_t R = \Delta_g R + 2|\Ric|_g^2 \geq \Delta_g R + \frac{2}{3} R^2
\end{equation}
which when combined with the standard maximum principle for scalar parabolic equations yields the useful monotonicity formula
\begin{equation}\label{rome}
 \partial_t R_\min \geq \frac{2}{3} R_\min^2
 \end{equation}
If $R_\min$ was initially positive (i.e. the manifold had strictly positive scalar curvature to begin with), this monotonicity formula implies finite time blowup of Ricci flow, a fact first observed by Hamilton \cite{ham1}.  When $R_\min$ is negative, the monotonicity formula is not quite as powerful, but is still useful as one approaches a blowup singularity, because the $R_\min$ is \emph{subcritical} with respect to scaling.  Indeed it is then easy to see that if one rescales a Ricci flow around its blowup time, then the scalar curvature becomes asymptotically non-negative.  This is a weak version of the \emph{Hamilton-Ivey pinching phenomenon}, of which more shall be said shortly.

A variant of $R_\min$ which is also monotone non-decreasing is the least eigenvalue of $-4\Delta_g + R$.  In fact Ricci flow can be viewed as a gradient flow for this quantity: see \cite[\S 1,2]{per1}, \cite[\S 8]{per2}, \cite[\S 5-9]{kleiner}, or \cite[\S 1.5]{cao}.  However we do not need this monotonicity formula here.  We also mention Hamilton's entropy functional $\int_M R \log R\ dg$, which is monotone in the case of positive curvature \cite{ham1}, although again this will not be needed here.

The Ricci tensor $\Ric$ and Riemann tensor $\Riem$ also obey equations similar to \eqref{r-eq}, provided one works in a suitably chosen orthonormal frame.  Indeed one has equations which essentially take the form
\begin{equation}\label{riem}
\partial_t \Riem = \Delta_g \Riem + O( \Riem^2 )
\end{equation}
in the orthonormal frame coordinates\footnote{The selection of such a frame modifies the dimensional analysis, basically by setting $G=1$, thus collapsing each of the columns of Table \ref{dimtable} to a single entry.}.  These equations can be combined with a tensor version of the maximum principle to obtain some uniform lower or upper bounds on the Riemann tensor, similar to those obtained for the Ricci scalar.  Some simple examples of this are the following:

\begin{itemize}

\item If the scalar curvature of a Ricci flow is non-negative at the initial time, then it remains non-negative at all later times.  (This is a special case of the monotonicity of $R_\min$.)

\item If the Ricci tensor of a Ricci flow is non-negative definite at the initial time, then it remains non-negative definite at all later times.  

\item If the sectional curvatures of a Ricci flow are all non-negative\footnote{An extremely minor remark: it would be better to replace the non-negativity of sectional curvatures by the non-negativity of the Riemann curvature tensor, as then this claim extends to all dimensions.  However for three dimensions the two statements are equivalent, due to the identification between two-forms and (scalar multiples of) planes in this setting.} at the initial time, then they remain non-negative at all later times.

\end{itemize}

See for instance \cite[\S 4]{morgan} or \cite[\S 2]{cao} for a complete discussion of these results, which essentially are due to Hamilton \cite{pinch}.  An important variant of the second fact is \emph{Hamilton's rounding theorem} \cite{ham1}, which asserts that if the Ricci curvature is initially strictly positive, then not only does the Ricci flow blow up in finite time, the metric also becomes asymptotically spherical as one approached the blowup time, thus establishing the Poincar\'e conjecture in the positive Ricci curvature case.

There are more quantitative versions of these results based on understanding the dynamics of the lowest eigenvalue of the Ricci tensor or Riemann tensor.  They roughly have the following form: if these tensors are \emph{nearly} non-negative (so that the lowest eigenvalue is not too negative, relative to the size of the tensor) at the initial time, then they remain nearly non-negative at later times.  In particular, even if the curvature is blowing up in magnitude, the negative portion of the curvature will not blow up as fast; this phenomenon is known as \emph{Hamilton-Ivey pinching towards positive curvature}.  One (informal) statement of this pinching is a lower bound on sectional curvatures $K = K(X,Y)$ which roughly takes the following form: if it is known that $K \geq -O(1)$ at time zero, then for later times $t$ (which are not too large) one has a bound of the form
\begin{equation}\label{kr}
 K \geq -o(1 + |R|)
 \end{equation}
where $o(1 + |R|)$ denotes some quantity which goes to infinity slower than $1 + |R|$ as $|R| \to \infty$ (note that monotonicity of $R_\min$ ensures that $R \geq -O(1)$).  This pinching 
also lets one deduce an upper bound
$$ |\Riem|_g, |\Ric|_g = O( 1 + |R| ),$$
because pinching prevents the scenario where large positive and large negative curvature cancel to create a small scalar curvature.  Thus the scalar curvature $R$ in fact controls in magnitude all the other curvature components, and is particularly effective at controlling the negative components of curvature.

This pinching statement was proven in \cite{pinch}.  Roughly speaking, it is proven by showing that for a carefully selected (and mildly time-dependent) choice of $o(1+|R|)$ expression in \eqref{kr}, the set of points obeying \eqref{kr} enjoys a certain convexity-in-time property; the hypotheses on initial conditions combined with a continuity argument then give the result.  As a very crude first approximation, the pinching phenomenon means that in high curvature regions, almost all the negative curvature has been eliminated, leaving only non-negatively curved geometries to consider.

One common application of maximum principles is to show that vanishing of a non-negative elliptic or parabolic object at one point implies vanishing at other points.  For instance, the classical maximum principle for harmonic functions implies that if a non-negative harmonic function on a domain vanishes at an interior point, then it vanishes everywhere.  Hamilton's tensor maximum principle implies a result of a similar flavour: if a three-dimensional Ricci flow with non-negative sectional curvatures has a vanishing sectional curvature at one point, then (locally, at least) the flow must split into the product of a two-dimensional Ricci flow, and a flat line $\R$.  This type of splitting result is used heavily in the classification of various types of limiting Ricci solution later on in the argument; see \cite[\S 4.4]{morgan}, \cite[\S 40]{kleiner} or \cite[\S 2.2]{cao}.

Another natural quantity which is controlled (though not monotone) is the total volume $\Vol(t) := \int_{M_t}\ dg(t)$.  Indeed, using the the formula for the rate of growth of the volume form under Ricci flow, followed by monotonicity of $R_\min$,
$$ \partial_t dg = - R dg \leq -R_{\min}(t) dg \leq -R_\min(0) dg,$$
we see from Gronwall's inequality that the volume of the manifold can grow at most exponentially in time, with the rate determined by initial minimal scalar curvature $R_\min(0)$.  Unfortunately the volume is a supercritical quantity and is not directly useful for blowup analysis, although we will need this volume control later in order to show that the blowup times are discrete (the key point being that each surgery removes a significant portion of volume from the manifold).

If the sectional curvatures are all non-negative, then much stronger monotonicity formulae are possible.  For instance, the above analysis shows that the volume is now monotone decreasing; some further monotonicity results of this nature are discussed in the next section.  Furthermore, there are stronger pointwise monotonicity formulae for the Ricci scalar $R$, arising from \emph{Hamilton's Harnack inequality} \cite{hamil-surf} (inspired by an earlier formula of Li and Yau \cite{liyau} for a different equation).  There are several formulations of this inequality.  One is the following: if the Ricci flow has non-negative bounded sectional curvature on $[0,t]$, then one has the pointwise monotonicity-type formula for the scalar curvature
$$ \partial_t R + \frac{1}{t} R + 2 \partial^\alpha R V_\alpha + 2 \Ric_{\alpha \beta} V^\alpha V^\beta \geq 0$$
for all vector fields $V$; this inequality can be obtained after some effort from many geometric identities and the maximum principle.  If $0 < t_1 < t_2$, one can then quickly deduce the pointwise inequality
$$ R(x_2, t_2) \geq \frac{t_1}{t_2} e^{-d_{M_{t_1}}(x_1,x_2)^2/2(t_2-t_1)} R(x_1,t_1)$$
which may help explain the terminology ``Harnack inequality''.  In particular, for ancient solutions (extending to arbitrarily negative times) with non-negative bounded sectional curvature, then $R(x,t)$ is a non-decreasing function of $t$, thus we have \emph{pointwise} monotonicity of the scalar curvature in this case.  See \cite[\S 7,9]{per1}, \cite[\S 4.5.2]{morgan}, \cite[Appendix F]{kleiner}, or \cite[\S 2.5]{cao} for details.

The above consequences of the maximum principle are quite robust.  For instance, they all extend (with minor modifications) to general dimensions.  Furthermore, while many of the stated results and arguments are initially only valid for compact manifolds, there are a number of ways to extend them to non-compact manifolds given a suitable control on the geometry at infinity.  We shall gloss over the technical details (usually involving suitable weight functions to localise arguments to near a compact set) in performing such extensions, but see \cite{shi} or \cite[\S 12.2]{morgan} for examples.

\section{Monotonicity formulae II: minimal surfaces and finite time extinction}\label{Extinct-sec}

In the previous section we saw how \emph{pointwise} quantities such as curvature enjoy monotonicity properties, largely thanks to the maximum principle.  In geometry, another fertile source of monotone (or at least interesting) quantities arises not from \emph{points} in the manifold, but from embedded \emph{curves} and \emph{surfaces}, particularly those which are chosen to minimise some geometric functional such as length, area, or energy.  It turns out that the Ricci flow tends to contract\footnote{Strictly speaking, these this contraction phenomenon only starts kicking in once the minimal Ricci curvature $R_\min$ stops being too negative, so this is only a conditional monotonicity formula.  However, by combining these formulae with the monotonicity formula \eqref{rome} for $R_\min$ we can still extract finite time blowup results.}  these sorts of quantities, which leads to a useful way to force finite time blowup (an important aspect of Theorem \ref{main}).  The intuition seems to be that the Ricci flow ``wants'' to compress all of the non-trivial higher-dimensional topology of the manifold (e.g. the homotopy groups $\pi_2(M)$ and $\pi_3(M)$) into nothingness.  This compression effect can be temporarily counteracted by large negative curvature, which can act to expand the metric, but the pinching phenomenon described in the previous section eventually weakens the negative curvature components of the geometry to the point where the compression becomes unstoppable\footnote{Actually, this is only true for the higher-dimensional topology.  For the first homotopy group $\pi_1(M)$, it is possible that the negative curvature never weakens to the point where compression occurs; this can be seen for instance by looking at a constant negative curvature manifold (which can have non-trivial $\pi_1(M)$ but not $\pi_2(M)$ or $\pi_3(M)$), in which Ricci flow expands the manifold indefinitely.  We thank John Morgan for this example.}.

In this section we discuss some monotone (or at least controlled) quantities associated to appropriately minimising one-dimensional manifolds (curves) and two-dimensional manifolds (surfaces).  The quantity $R_\min$ in the preceding section can be thought of as a zero-dimensional case of these scheme, while the total volume $\int_M\ dg$ is a three-dimensional case\footnote{Perelman's monotone quantities, discussed in the next section, have an interpretation which resembles an ``infinite-dimensional'' quantity, though it is unclear to me how they exactly fit into the energy-minimising philosophy.}.

To give an example of this compression effect for one-dimensional minimal manifolds (i.e. geodesics), let us observe the following simple but useful fact:

\begin{proposition}[Non-negative curvature implies distance contraction]\label{nonneg}  Suppose that a Ricci flow $t \mapsto M_t$ has non-negative sectional curvature.  Then for fixed $x,y \in M$, the distance function $d_t(x,y)$ is a non-increasing function of $t$.
\end{proposition}

\begin{proof}  The first variation formula for metrics, combined with the Ricci flow equation, yields
$$ \partial_t \frac{1}{2} d_t(x,y)^2 = - \int_0^1 \Ric_{\alpha \beta}(\gamma(s)) \dot \gamma^\alpha(s) \dot \gamma^\beta(s) ds$$
for some energy-minimising curve $\gamma: [0,1] \to M$ from $x$ to $y$.  The claim follows.
\end{proof}

Clearly one can extend and refine this type of estimate; see for instance \cite[\S 8]{per1}, \cite[\S 3.6]{morgan}, \cite[\S 3.4]{cao}, or \cite[\S 26]{kleiner} for versions of the above compression phenomenon of relevance to Perelman's argument.

Using Proposition \ref{nonneg} one could then create a monotonic decreasing functional in the non-negatively curved case, namely the minimal length of a non-contractible closed loop.  Unfortunately, this quantity is only non-zero when the homotopy group $\pi_1(M)$ is non-trivial - and for the purposes of proving the Poincar\'e conjecture, this is certainly not the case!  However, let us note that this result will yield a finite time extinction result for non-simply connected manifolds with strictly positive curvature (though this was already known to Hamilton by other means).
Thus we must move to surfaces (such as embedded two-spheres or two-tori) to get some interesting geometric quantities.  Experience has shown that surfaces which minimise area or energy\footnote{Area and energy of surfaces are closely related, but not quite the same.  A surface $\beta: \R^2 \to M$ has area $\int_{\R^2} |\partial_x \beta \wedge \partial_y \beta|_g\ dx dy$ and energy
$\frac{1}{2} \int_{\R^2} |\partial_x \beta|_g^2 + |\partial_y \beta|_g^2\ dx dy$.  By Cauchy-Schwarz, the area is thus less than the energy, with equality when $\beta$ is conformal.  In practice, surfaces such as spheres can usually be made conformal after reparameterisation, and so minimising the area and minimising the energy usually end up amounting to the same thing.}
are likely be of great use in geometry; to pick just one example, witness the use of minimal surfaces in Gromov's ground-breaking paper \cite{gromov} on symplectic non-squeezing.  Actually, it turns out that minimising over surfaces only is able to exploit $\pi_2(M)$ information; it lets one reduce to the case where $\pi_2(M)$ is trivial, but this turns out to not be enough.  One may think then to move to embedded three-spheres to exploit $\pi_3(M)$ information, but this turns out to be difficult to do directly.  Fortunately, we can view embedded three-spheres as a loop of two-spheres, or as a two-sphere of loops, and by taking an appropriate min-max over the two parameters one obtains a useful monotone quantity that allows one to exploit $\pi_3(M)$ information, which turns out to be sufficient to force finite time extinction.  The loop-of-two-spheres approach was pursued by Colding and Minicozzi \cite{cold}, whereas the two-sphere-of-loops approach was pursued by Perelman \cite{per3}.  We shall discuss the former which is slightly simpler conceptually, even though for technical reasons the proof of Poincar\'e's conjecture turns out to be more easily pursued using the latter; see \cite[\S 18]{morgan} for details\footnote{The papers \cite{cao}, \cite{kleiner} do not address the finite time extinction issue as they are following Perelman's proof of the geometrisation conjecture, which does not require such extinction.}.

We remark that the idea of using minimal surface areas as controlled quantities under Ricci flow originates from Hamilton \cite{pinch}.

Let us first discuss the $\pi_2(M)$-based quantities.  Suppose that $\pi_2(M_t)$ is non-trivial.  Let $\beta: S^2 \to M_t$ be any immersed sphere not homotopic to a point.  Each such sphere has an energy $E(\beta,t) := \frac{1}{2} \int_{S^2} |d\beta|_{g_t}^2$ using the metric $g_t$ at time $t$; this is at least as large as the area of $\beta$.  Define $W_2(t)$ to be the infimum of $E(\beta,t)$ over all such $\beta$.  It turns out (from standard Sacks-Uhlenbeck minimal surface theory \cite{sacks}) that this infimum is actually attained (at which point the energy is the same as the area).  A computation (using the local structure of minimal surfaces, and the Gauss-Bonnet formula) lets one obtain the differential inequality
$$ \partial_t W_2(t) \leq -4\pi - \frac{1}{2} R_\min(t) W_2(t)$$
(where the derivative is interpreted in a suitably weak sense) whenever the Ricci flow remains smooth.  Combining this with \eqref{rome} and some elementary Gronwall-type analysis one can show that $W_2(t)$ becomes negative in finite time, which is absurd.  This forces blowup in finite time whenever $\pi_2(M)$ is non-trivial.  It turns out that a similar argument can also be extended for Ricci flow with surgery.  The final conclusion is that after finite time, the Ricci flow with surgery is either completely extinct, or at every time each component of $M$ trivial $\pi_2$.

A variant of this argument works for $\pi_3(M)$.  If $\pi_3(M)$ is non-trivial, then there exists a loop of two-spheres $\gamma: [0,1] \times S^2 \to M$ with $\gamma(0) = \gamma(1)$ equal to a point, which cannot be contracted to a point.  For any such $\gamma$, define $E(\gamma,t) := \sup_{0 \leq s \leq 1} E(\gamma(s),t)$ be the largest energy obtained by a two-sphere on the loop, and then let $W_3(t)$ be the infimum of $E(\gamma,t)$ over all such $\gamma$.  A variant of the above argument then gives\footnote{This is assuming some analytical facts concerning near-extremisers of $W_3(t)$ which are widely expected to be true, but have not yet been explicitly proven.
More precisely, one needs a bubbling-type concentration compactness property for minimising sequences of loops of spheres, which is likely to be true (based on analogy with similar results for simpler variational problems).  Nevertheless it is believed that this property will be established soon.  We thank Gang Tian and John Morgan for clarifying these points.}
$$ \partial_t W_3(t) \leq -4\pi - \frac{1}{2} R_\min(t) W_3(t)$$
so once again we can demonstrate finite time blowup as long as $\pi_3(M)$ is non-trivial.  Again, it is possible (given sufficient care) to extend these arguments to the case of Ricci flow with surgery.

There is the issue of $\pi_3(M)$ being non-trivial.  Fortunately some homology theory settles this fairly readily, at least in the simply connected case.  If $M$ is simply connected, then $\pi_1(M)$ is trivial, which implies (by the Hurewicz theorem) that the first homology group $H_1(M)$ is also trivial.  Poincar\'e duality then gives that $H_2(M)$ is also trivial, whereas the three-dimensionality\footnote{Note that simply connected manifolds are necessarily orientable.} of $M$ gives that $H_3(M) = \Z$.  In short, $M$ has the same homology groups as the three-sphere, and thus by Whitehead's theorem it is therefore homotopy equivalent to a sphere\footnote{Note that this is a substantially weaker statement than being \emph{homeomorphic} to a sphere; indeed, since non-trivial homology spheres are known to exist, the Poincar\'e conjecture cannot be proven by homology theory alone!}.  In particular\footnote{Here we are using the fact that the surgery procedure cannot destroy the simply connected nature of $M_t$ (in fact the fundamental groups of the post-surgery components can easily be shown to be subgroups of the pre-surgery fundamental groups, which are trivial in the simply connected case).} it quickly follows that $\pi_3(M)$ is non-trivial (indeed, it is isomorphic to $\Z$).  
The more general case when $\pi_1(M)$ is a free product of finite groups and infinite cyclic groups can be handled by more sophisticated topological and homological arguments of the above type; see \cite[\S 18]{morgan} for details.

As mentioned earlier, the Colding-Minicozzi approach requires a technical analytical property which is expected to be true but whose proof is not currently in the literature.  There is an alternate approach of Perelman \cite{per3} for which the corresponding analytical property is easier to verify (because one works with area-minimising surfaces instead of arbitrary surfaces).  Instead of viewing embedded three-spheres $\gamma$ in $M$ as loops of two-spheres, one instead views them as two-spheres of loops (with some natural boundary conditions).  Each loop spans a minimal disk (because of simple connectedness), which of course has an area (or energy).  Taking the maximal area obtained in this way along a two-sphere of loops, one obtains an energy $\tilde E(\gamma,t)$ analogous to the energy considered earlier; maximising over $t$, and then minimising over all non-contractible $\gamma$ one obtains an analogue $\tilde W_3(t)$ to $W_3(t)$, which ends up obeying a very similar differential inequality\footnote{A significant technical issue arises, which is that to demonstrate the above inequality it becomes necessary to evolve the loops involved by curve shortening flow.  This flow can become singular even when the Ricci flow is smooth, requiring a (standard) lift of the problem to higher dimensions to resolve the singularities.  See \cite{per3}, \cite[\S 18.4-18.6]{morgan} for details.} to $W_3(t)$.  
For full details see \cite[\S 18]{morgan}.

In summary, modulo some important technical details (such as how to deal with surgery), we have used the compression of minimal-area surfaces phenomenon to obtain the finite time extinction property for Ricci flow-with-surgery.  In the bounded amount of time
remaining to the flow, we need to control the singularities (or near-singularities) well enough that 
we can correctly define a surgery procedure which (a) does not disrupt other aspects of the argument, such as the compression phenomenon, and (b) modifies the topology of the manifold in a controllable manner.  These (highly non-trivial) tasks will be the focus of the remaining sections of this note.

\section{Monotonicity formulae III: Perelman's reduced volume}\label{perelcrit}

One of Perelman's major advances in the subject was to introduce for the first time in \cite[\S 3,7]{per1} several useful \emph{critical} and \emph{coercive} monotone quantities for Ricci flow, most notably\footnote{Perelman also introduces a number of other interesting critical monotone quantities, such as a modified entropy functional ${\mathcal W}$.  While these quantities are also very useful, particularly for the geometrisation conjecture, they are not strictly necessary for Theorem \ref{main} and will not be discussed here.}
 a \emph{reduced volume} $\tilde V$ which is a ``renormalized'' version of the total volume.  This quantity is not only scale-invariant, but is also monotone under Ricci flow with no assumptions on the  initial signature of the Ricci curvature.  This is an absolutely essential breakthrough if one is to analyse behaviour near a singularity, or more generally near a high curvature region of the manifold.

Now we discuss the reduced volume in more depth\footnote{Details on the reduced volume can be found in \cite[\S 7]{per1}, \cite[\S 6]{morgan}, \cite[\S 14-23]{kleiner}, or \cite[\S 3]{cao}.}.  This quantity is consistent with the standard ``heat ball'' method for generating monotonicity formulae for geometric parabolic equations.  Note that for a solution to the free heat equation $u_t = \Delta u$ in a Euclidean space\footnote{Actually, the monotonicity formulae given here are valid (and scale-invariant) in all spatial dimensions, thus raising the possibility that Perelman's methods should extend to analyse Ricci flow on manifolds in every dimension; see for instance \cite{chenzhu} for some preliminary progress in this direction.  However, at present several other aspects of Perelman's work is restricted to three spatial dimensions.} $\R^3$, and any point $(t_0,x_0)$ in spacetime, the quantity
\begin{equation}\label{heatball}
 \frac{1}{(4\pi\tau)^{3/2}} \int_{\R^3} e^{-\dist_{\R^3}(x,x_0)^2/4\tau} u(t_0-\tau,x)\ dx
 \end{equation}
is constant in time for $\tau > 0$; this is clear from the fundamental solution of the heat equation.  It turns out that when one replaces the free heat equations by other parabolic equations, that similar quantities to \eqref{heatball} can often turn out to be monotone in time rather than constant.  As a gross caricature, one can think of \eqref{heatball} as essentially measuring an average of $u(t_0-\tau)$ on a ball centred at $x_0$ of radius $O( \sqrt{\tau} )$.  Observe also that if we normalise $(t_0,x_0)$ to the origin, then \eqref{heatball} is also invariant under the scaling $u(t,x) \mapsto u(\frac{t}{\lambda^2}, \frac{x}{\lambda})$, which is the analogue to \eqref{local-scale}.

Let $(t_0,x_0)$ be a point in spacetime of a solution $M$ to Ricci flow in $d$ spatial dimensions, thus $x_0$ lies 
on the $3$-dimensional Riemmanian manifold $(M_{t_0}, g_{t_0})$.  Based on the above discussion, one might naively expect the expression
\begin{equation}\label{av-heat}
 \frac{1}{(4\pi\tau)^{3/2}} \int_{M_{t_0-\tau}} e^{-\dist_{M_{t_0-\tau}}(x,x_0)^2/4\tau}\ dg_{t_0-\tau}(x)
\end{equation}
 to have some monotonicity in time; this is basically measuring the volume of the ball of radius $O(\sqrt{\tau})$ centred at $x_0$ at time $t_0-\tau$, divided by the volume of the corresponding Euclidean ball.
Thus this quantity will be expected to equal one when $M$ is flat, be less than one when $M$ has positive curvature,
and greater than one when $M$ has negative curvature.  Note also when we normalise $(t_0,x_0)$ to be the origin, then
the quantity $\dist_{M_{t_0-\tau}}(x,x_0)^2/4\tau$ is dimensionless (scale-invariant), as is the quantity \eqref{av-heat}.

Unfortunately the above quantity is not quite monotone.  Geometrically, this can be explained by the fact that $x_0$ ``lives'' on the manifold $M_{t_0}$ at time $t_0$, and really shouldn't be used to measure distances on the manifold
$M_{t_0-\tau}$ at time $t_0-\tau$.  Remarkably, however, it turns out that a slight modification of the above functional \emph{is} monotone non-increasing (as a function of $\tau$) for smooth Ricci flow as $\tau \to 0^+$, namely the \emph{reduced volume}\footnote{This definition differs from the one by Perelman by a harmless factor of $(4\pi)^{3/2}$; I have chosen this normalisation to give the flat metric a reduced volume of $1$.}
$$
\tilde V(\tau)  = \tilde V_{t_0,x_0}(\tau):= \frac{1}{(4\pi\tau)^{3/2}} \int_{M_{t_0-\tau}} e^{-l_{t_0,x_0}(t_0-\tau, x)}\ dg_{t_0-\tau}(x)$$
where the \emph{reduced length} $l_{t_0,x_0}(t_0-\tau, x)$ is a variant of the naive Gaussian distance
$\dist_{M_{t_0-\tau}}(x,x_0)^2/4\tau$, defined as follows.  Observe from the definition of geodesic distance and elementary calculus of variations that
$$\frac{\dist_{M_{t_0-\tau}}(x,x_0)^2}{4\tau} = \inf_\gamma \frac{1}{2\sqrt{\tau}} 
\int_0^\tau |\dot \gamma(\tau')|_{g_{t_0-\tau}}^2\ \sqrt{\tau'} d\tau'$$
where $\gamma: [0,\tau] \to M_{t_0-\tau}$ runs over all smooth curves with $\gamma(0) = x_0$ and $\gamma(\tau) = x$.
(Indeed, the minimiser is attained when $\gamma$ traverses the geodesic from $x_0$ to $x$, which at time $\tau'$ is at distance $\sqrt{\frac{\tau'}{\tau}} \dist_{M_{t_0-\tau}}(x,x_0)$ from $x_0$.)  As mentioned before, it is not geometrically natural for $x_0$ to be placed in $M_{t_0-\tau}$ when it really should lie in $M_{t_0}$; similarly, it is not natural to $\gamma$ to be completely contained in $M_{t_0-\tau}$, and instead it should be a curve from $(t_0,x_0)$ to $(t_0-\tau,x)$.  This (at least partially) motivates Perelman's definition of reduced length, namely
$$ l_{t_0,x_0}(t_0-\tau,x) := \frac{1}{2\sqrt{\tau}} {\mathcal L}_{t_0,x_0}(t_0-\tau,x)$$
and ${\mathcal L}$ is the \emph{length}
\begin{equation}\label{length-func}
 {\mathcal L}_{t_0,x_0}(t_0-\tau,x) = \inf_\gamma \int_0^\tau [|\dot \gamma(\tau')|_{g_{t_0-\tau'
}}^2 + R(t_0-\tau', \gamma(\tau'))]\ \sqrt{\tau'} d\tau'
\end{equation}
where $\gamma: [0,\tau] \in M$ ranges over all smooth curves with $\gamma(\tau') \in M_{t_0 - \tau'}$ for $0 \leq \tau' \leq \tau$, with $\gamma(0) = x_0 = (t_0,x_0)$ and $\gamma(\tau) = x = (t_0-\tau, x)$, and $R$ is the Ricci scalar at the point $(t_0-\tau',\gamma(\tau') \in M_{t_0-\tau'}$.  
This type of ``parabolic length'' had been considered earlier in a different context than 
Ricci flow by Li-Yau \cite{liyau}.  The presence of the Ricci scalar appears to be some sort of correction term, similar in spirit to the Ricci curvature terms arising in other geometric formulae such as the second variation formula for geodesic distance, the Bochner-Wietzenb\"ock formula for the Laplacian of energy densities, or the presence of the Ricci scalar in the conformal Laplacian.  See \cite[\S 7]{per1}, \cite[\S 2.5]{cao} for some motivation of this length functional arising naturally from a Li-Yau type inequality for the backwards heat equation.  

One can easily verify that $l_{t_0,x_0}(t_0-\tau,x)$ is dimensionless, and hence the reduced volume $\tilde V(\tau)$ is also scale-invariant.  Thus the monotonicity of this quantity is thus a \emph{critical} control on the geometry as one approaches a blowup point $(t_0,x_0)$; very informally speaking, it asserts that as one approaches such a singularity, the net curvature tends to increase, thus there is a trend from negative curvature to positive.  In some ways this can be viewed as a local version of the monotonicity formula for $R_\min$.  

The proof of this monotonicity is remarkably delicate; see \cite[\S 7]{per1}, \cite[\S 22]{kleiner}, \cite[\S 6.7]{morgan}, or \cite[\S 3.2]{cao} for the details.  In principle, one simply needs to differentiate the reduced volume in time, use variation formulae to compute various derivatives of the reduced length, and integrate by parts, using many identities from Riemannian geometry, until one visibly sees the sign of the derivative.  In practice, of course, such a computation is very lengthy and seems highly miraculous in the absence of geometric intuition.  Since I do not have this intuition, I cannot really shed any light as to why this monotonicity formula works (other than by analogy with simpler equations).  It does however appear to be instructive to show that the monotone quantity is in fact constant in the case of a gradient shrinking soliton.  In \cite[\S 6]{per1}, \cite[\S 3.1]{cao} an interpretation of reduced volume in terms of infinite dimensional Riemannian geometry was given, relating this quantity to a matrix Harnack expression considered by Hamilton \cite{hamil-harnack}.  A related (but not identical) geometric interpretation for this expression also appears in \cite{cc}.

Comparing the reduced volume against the informal interpretation for \eqref{av-heat} leads to the heuristic
\begin{equation}\label{tv}
 \tilde V_{t_0,x_0}(\tau) \approx \tau^{-3/2} \Vol_{M_{t_0-\tau}}( B(x_0, \tau^{1/2}) ).
\end{equation}
This heuristic turns out to be essentially accurate in the case where the ball $B(x_0, \tau^{1/2})$ has
bounded normalised curvature on the time interval $[t_0-\tau,t_0]$, in the sense that $|\Riem|_g = O( \tau^{-1} )$ (this is the natural amount of curvature predicted by dimensional analysis, and is also consistent with \eqref{riem}).  Verifying this is not trivial, but roughly speaking proceeds by exploiting the curvature control to ensure the reduced length $l_{t_0,x_0}(t_0-\tau,x)$ stays bounded on $B(x_0,\tau^{1/2})$,
which in turn leads to an estimate similar to \eqref{tv}.

A more refined analysis of the proof of the monotonicity formula also reveals the precise solutions for which equality occurs, i.e. the reduced volume is constant rather than monotone.  Indeed, this occurs precisely when the solution is a \emph{gradient shrinking soliton}, as in \eqref{gss} (if we normalise $t_0=0$ and $t=-1$).  This strongly suggests (but does not formally imply) that as one rescales around a singularity $(t_0,x_0)$ of a Ricci flow, that the flow should increasingly resemble a gradient shrinking soliton (possibly after passing to a subsequence of rescalings).  A result of this form will be essential later in the argument.  Details can be found\footnote{The argument in \cite[\S 7]{per1} is slightly incorrect, however can be fixed in a number of ways; see the other three references cited here.} in \cite[\S 7]{per1}, \cite[\S 9.2]{morgan}, \cite[\S 36]{kleiner}, \cite[\S 3]{cao}.

We have normalised the reduced volume to equal $1$ when the manifold is flat.  If the Ricci flow has bounded curvature and is smooth at the final point $(t_0,x_0)$, then its rescalings become flat, and so the reduced volume is asymptotically $1$ as $t \to t_0$.  By monotonicity this means that the reduced volume is less than or equal to $1$ at all times $t < t_0$.  By analysing when equality occurs in the monotonicity formula as in the previous paragraph, one can in fact conclude that the reduced volume will be \emph{strictly} less than $1$ unless the entire flow was flat.  See \cite[\S 9.2]{morgan} for details.

\section{Blowup analysis I: a no-local-collapsing theorem}\label{nlc}

Armed with a critical monotone quantity centred at a spacetime point of our choosing, we can now hope to obtain control on the manifolds $M_t$ of a Ricci flow\footnote{In order to fully prove Theorem \ref{main}, one has to consider generalised Ricci flows, in which various surgeries take place between $t=0$ and $t=t_*$.  This makes the analysis significantly more complicated and non-trivial, but for sake of this high-level overview we shall focus on the much simpler case of smooth Ricci flow (without surgeries).} as one approaches a blowup time, $t \to t_*$.

This will be done in several stages.  The first major stage, which we cover in this section, involves the relationship between the ``local'' geometry of the manifold $M_t$, which is given by pointwise quantities such as the Riemann curvature tensor $\Riem$ (and its components such as the Ricci tensor $\Ric$ and Ricci scalar $R$), and more ``global'' geometry of the manifold, such as the injectivity radius or volume of balls $B(x,r)$ at a macroscopic radius $r$.  The local quantities will be controllable directly from equations arising from the Ricci flow such as \eqref{riem}, whereas it is the global quantities which will be needed to control the asymptotic profile of blowup solutions and perform surgery properly.
Given an arbitrary smooth Riemannian manifold there are some general partial relationships between the local and global quantities.  Let us first work at the macroscopic radius scale $r=1$.  Just to give the flavour of things we shall use informal $O()$ notation (thus suppressing the exact dependence of various implicit constants on other implicit constants, which is an important but technical detail) and not give the sharpest results:

\begin{itemize}

\item (Comparison theorem) If one has bounded curvature on a unit ball $B(x_0,1)$, thus $|\Riem|_g = O(1)$ on this ball, then for $0 \leq r \leq 1$, the volume of $B(x_0,r)$, divided by $r^3$, is essentially non-increasing in $r$.  In particular this implies that $B(x_0,1)$ has volume $O(1)$.

\item (Reverse comparison theorem)  If one has bounded curvature on a unit ball $B(x_0,1)$, then the volume of the ball is $\sim 1$ if and only if the injectivity radius\footnote{Roughly speaking, the injectivity radius at a point $x$ is the radius of the largest ball $B(x,r)$ for which length-minimising geodesics are unique.} is $\gtrsim 1$.

\end{itemize}

These results are facts from Bishop-Cheeger-Gromov comparison theory and rely primarily on an analysis of the exponential map from $x_0$.  The curvature bounds ensure that this map has controlled derivatives, leaving the injectivity of the map as only remaining hypothesis required to estimate the volume accurately.  For the more precise versions needed for Perelman's argument, see \cite[\S 6.2]{petersen}, \cite[\S 1.6-7]{morgan}, or \cite[\S 4.7]{cgt}.

Note that on a generic smooth Riemannian manifold it is possible to have \emph{collapsing}, in which the curvature is bounded but the volume of the ball is small (and the injectivity radius is small too, by the reverse comparison theorem).  A simple example is given by an extremely small flat torus or flat\footnote{This is of course flatness in the intrinsic sense of Riemann, as opposed to any sort of extrinsic sense; indeed, we are not actually embedding our manifolds into an ambient Riemannian space with which to measure extrinsic curvatures such as second fundamental forms.} cylinder, which has zero curvature but for which the unit ball has unexpectedly small volume.  Many more exotic examples can of course be constructed.  In contrast, in the non-collapsed case, a famous finiteness theorem of Cheeger \cite{cheg} limits the number of topologically distinct such examples available.  Thus we see that the non-collapsed property is very powerful in constraining the behaviour of manifolds.  As the small torus and cylinder examples show, collapsing is associated with the presence of ``little loops'' (i.e. short closed geodesics).

Of course one can rescale the above results to work at any radius scale $r$, to obtain:

\begin{itemize}

\item (Comparison theorem) If one has bounded normalised curvature on a ball $B(x_0,r)$, thus $|\Riem|_g = O(r^{-2})$ on this ball, then for $0 \leq r' \leq r$, the volume of $B(x_0,r')$, divided by $(r')^3$, is essentially non-increasing in $r'$. In particular $B(x_0,r)$ has volume $O(r^3)$.

\item (Reverse comparison theorem)  If one has bounded normalised curvature on a ball $B(x_0,r)$, then the volume of the ball is $\sim r^3$ if and only if the injectivity radius is $\gtrsim r$.

\end{itemize}

One consequence of the monotonicity aspect of the comparison theorem is that if there is no collapsing at a scale $r$, then there is also no collapsing at any smaller scale\footnote{It may be more intuitive to view this assertion in the contrapositive; if there is some breakdown of injectivity at a small scale caused by little loops, then this will replicate itself at larger scales also, since little loops can be concatenated to form big loops.}.

A remarkable phenomenon concerning Ricci flow is that as one approaches a blowup time $t \to t_*$, collapsing cannot actually occur at scales $r = O( \sqrt{t_*-t} )$ which are finer than the natural length scale $\sqrt{t_*-t}$ associated to the near-blowup time $t$ (this is the scale given by dimensional analysis).  This no-local-collapsing theorem was first established by Hamilton \cite{hamil-sing} in the case of non-negative sectional curvature, and then in general by Perelman \cite{per1}.  We state here a fairly weak and imprecise version of this theorem (stronger versions are certainly available, see \cite[\S 4,7,8]{per1}, \cite[\S 7]{per2}, \cite[\S 8]{morgan}, \cite[\S 12, 25, 27, 78]{kleiner}, \cite[\S 3.3, 3.4, 7.6]{cao}):

\begin{theorem}[No local collapsing]\cite{per1}  Let $t \mapsto M_t$ be a smooth Ricci flow for $0 \leq t \leq t_0 = O(1)$, let $x_0 \in M$, and suppose we have noncollapsed geometry at time zero in the sense that the curvature is bounded by $O(1)$ and the volume of any unit ball is comparable to $1$ at time zero.  Then we have local noncollapsing at time $t_0$ in the following sense: given any $0 < r < \sqrt{t_0}$, if one has bounded normalised curvature on the parabolic cylinder
$[t_0-r^2,t_0] \times B(x_0,r)$, then the volume of the ball $B(x_0,r)$ at time $t_0$ is  comparable to $r^3$ (or equivalently, the injectivity radius of this ball at time $t_0$ is comparable or larger to $r$).
\end{theorem}

\begin{proof}(Sketch)  A rescaling argument using the comparison theorem lets us reduce to the case $t_0 \geq 1$ and $r < \sqrt{t_0}/2$. 
Because curvature is bounded at time zero, one can use \eqref{riem} to show that it will stay bounded for short times $0 \leq t \ll 1$.  This will give us enough control on the geometry to obtain a lower
bound $\tilde V_{t_0,x_0}(0) \gtrsim 1$ for the reduced volume at time $0$.  By
monotonicity this implies $\tilde V_{t_0,x_0}(t_0-r^2) \gtrsim 1$.  Using \eqref{tv} (pretending for sake of argument that this heuristic is accurate; the actual argument is more technical) and the comparison theorem we obtain the claim.
\end{proof}

There are many sharper versions of this no local collapsing theorem in Perelman's work, and also in
subsequent papers.  In particular, in order to fully complete the proof of Poincar\'e's conjecture one needs to obtain a no local collapsing theorem which is true for Ricci flow \emph{with surgery}, and not just smooth Ricci flow; this turns out to be significantly more complicated and delicate.  But the flavour of all these results remains unchanged: if there is no collapsing at initial times, then at later times one still has no collapsing as long as one is working at local scales (less than the natural length scale dictated by the time left until blowup) and the curvature is
controlled in a scale-invariant fashion.

\section{Blowup analysis II: asymptotic profiles}

One major reason for wanting a no-collapsing theorem is that it enables one to exploit \emph{Hamilton's compactness theorem} \cite{hamil-compact}: a sequence of Ricci flows\footnote{Strictly speaking, each Ricci flow needs to also come with an ``origin'' or ``base point'' to allow one to properly define convergence; it suffices to have uniform bounds on the geometry on an expanding sequence of balls around this origin, with radii going to infinity.} which are uniformly non-collapsed with uniformly bounded normalised curvature at a certain scale, have a convergent subsequence\footnote{We will be vague about exactly what ``convergence'' of a sequence of Ricci flows means, but (due to the parabolic smoothing effects of Ricci flow) one can largely rely on smooth convergence on compacta in the geometric sense, after picking a base point as mentioned earlier.}, though the limit might be a non-compact manifold.  This is basically because the non-collapsing yields injectivity radius lower bounds which can be used to give a common coordinate chart to all the manifolds (cf. \cite{cheeger}); also some parabolic regularity estimates of Shi (see \cite{shi-thesis}, \cite[Appendix B]{lutian}, \cite{chow}) give control on higher derivatives of curvature (essentially obtained by differentiating equations such as \eqref{riem} and then using the energy estimate and Sobolev embedding) which are sufficiently strong to give compactness (cf. the Arzela-Ascoli theorem).  Such compactness, when combined with standard rescaling arguments, allows one (at least in principle) to analyse the asymptotic behaviour of singularities, which is an important milestone in the analysis of any nonlinear PDE and in particular is essential for defining the surgery procedure for Ricci flow.

Let us first see how this would work in a best-case scenario.  Suppose we have a Ricci flow $t \mapsto M_t$ on $[0,T_*)$ which blows up at a finite time $T_*$.  Since the scalar curvature controls all the other curvatures due to pinching, we thus see that there is a sequence $t_n$ of times going to $T_*$ and $x_n \in M$ such that the scalar curvatures $R_n := R(t_n,x_n)$ go to infinity.  We make the optimistic assumption\footnote{Of course, by choosing $x_n$ to be near the point where the maximum is attained, one can obtain this assumption automatically.  But to understand surgery, we will need eventually to place $x_n$ next to where a \emph{specific} singularity occurs; the problem is that as the curvature near that singularity is becoming large, it may be that somewhere else there is another singularity where the curvature is becoming even larger still.  Fortunately this can be dealt with because the effect of each high curvature region ends up being localised, though to make the entire argument rigorous requires a slightly tricky downward induction on curvature; in other words, one needs to control the very high curvature regions first before attending to the moderately high curvature reasons.  See the next few sections for more discussion.} that at each time $t_n$, the curvature $R_n$ is close to the largest curvature that has yet been encountered, thus $R(t,x) \lesssim R_n$ for all $0 \leq t \leq t_n$.
  One consequence of this is that we automatically enjoy bounded normalised curvature at scales $R_n^{-1/2}$ and below and all times up to $t_n$.  By the no local collapsing theorem we can then also conclude that there is no collapsing at these scales or below up to time $t_n$.

We can then create rescaled flows $t \mapsto M^{(n)}_t$ on the interval $[-t_n R_n, (T_*-t_n)R_n)$ by setting
$$ (M^{(n)}_t, g^{(n)}_t) := (M_{t_n + t/R_n}, R_n g^{(n)}_{t_n + t/R_n}).$$
These rescaled flows are also Ricci flows, and their backwards time of existence extends to $-\infty$ as $n \to \infty$.  Also, because of our rescaling, we now have bounded normalised curvature and no collapsing at unit scales.  Thus we may apply Hamilton's compactness theorem (declaring $x_n$ as the origin of each $M^{(n)}_t$) and extract a convergent subsequence whose limit is an \emph{ancient} Ricci solution, i.e. one which has existed on the entire interval $(-\infty,0]$.  Furthermore this ancient solution is also non-collapsed at unit scales and below; and in fact if the $x_n$ are all piled on top of each other, one can compare the rescaled solutions to each other to obtain non-collapsing at \emph{all} scales for which the normalised curvature is under control.  We caution that even if the original manifold $M$ is compact, it is perfectly possible for the limiting ancient solution to be non-compact.  Furthermore, the Hamilton-Ivey pinching phenomenon ensures that the limiting ancient solution will have non-negative sectional curvature everywhere.

Roughly speaking, what this analysis suggests is that every high-curvature region of a Ricci flow should resemble a rescaled ancient non-collapsing solution to Ricci flow.  If we knew this, and we also knew enough about the geometry of ancient non-collapsing solutions, then we could hope to perform a controllable surgery on high-curvature regions in order to reduce their curvature while controlling the change in topology (and not disrupting the finite time extinction argument which will arise later).  However, the above argument required that at any given time, the high-curvature region one is studying has higher curvature than any other region previously encountered.  This is unfortunately not always the case (consider for instance two simultaneous singularities, one with a faster curvature blowup rate than the other).  Fortunately, it is possible to \emph{induct} on the curvature, starting with the highest curvature region and controlling that, then moving onwards to the next highest curvature region, and so forth.  In order to make this induction work, we cannot quite decouple the two major steps in the classification of curvature regions, namely the approximation by ancient non-collapsing solutions and the understanding of the geometry of these ancient solutions; they have to be applied alternately because we need enough control on the very high curvature regions to guarantee that they will not disrupt the analysis of the less high curvature regions.  For this reason it makes sense to talk about the ancient solutions first, before explaining why they model singularities.

Interestingly, we will eventually have to perform a secondary blowup analysis inside the first, in that the ancient solutions themselves are not the irreducible components of asymptotic behaviour, and contain inside them an even more special asymptotic class of solution, namely the \emph{gradient-shrinking solitons}.  The latter were completely classified by Perelman \cite[\S 11]{per1} (assuming, as we may, non-negative sectional curvature and non-collapsedness), which leads to a qualitative description of the ancient solutions, which in turn will lead to a qualitative description of the high curvature regions.  It seems that this two-step procedure is necessary to fully understand high curvature regions.  Consider for instance a high-curvature region which resembles a rescaled cigar.  Then while this cigar itself contains an asymptotically cylindrical piece at infinity (and the cylinder is a gradient-shrinking soliton), the compact end of this cigar does not have any model as a gradient-shrinking soliton.

\section{Blowup analysis III: the geometry of ancient non-collapsed solutions}

As discussed above, we shall skip for now the analysis of why high-curvature regions resemble
ancient non-collapsed solutions, and instead analyse these solutions directly.  Following Perelman\cite[\S 11]{per1}, \cite[\S 1]{per2} (and \cite[\S 9]{morgan}, \cite[\S 6]{cao}, \cite[\S 37]{kleiner}), let us define a\footnote{The name derives from a parameter $\kappa$ used to quantify the exact nature of the non-collapsing, but we are using the imprecise $\gtrsim$ notation to suppress this parameter.} \emph{$\kappa$-solution} to be a non-collapsed ancient Ricci flow of non-negative sectional curvature\footnote{Strictly speaking, one also has to hypothesise that the scalar curvature is pointwise non-decreasing, and that the curvature is bounded.  But the former follows from Hamilton's Harnack inequality of Li-Yau type discussed in earlier sections, while the latter will be established shortly.}, which is not flat; these will eventually be the solutions used to model the high-curvature regions.  The non-flatness will guarantee (in conjunction with some maximum principle arguments) that the scalar curvature $R$ is \emph{strictly} positive.
There are three model classes of examples of $\kappa$-solutions (and some non-examples) to keep in mind:

\begin{itemize}

\item (The sphere) For all $-\infty < t < 0$, let $M_t$ be the standard sphere $S^3$ with metric $g_t = (1-t) g_0$, where $g_0$ is the standard metric (rescaled to have constant sectional curvature $1$).  This is an ancient solution (in fact, a gradient shrinking soliton) starting with an infinitely large round sphere shrinking down to the unit sphere; if one continued the Ricci flow to positive times, then the sphere would shrink to a point at time one.  More generally one can replace the sphere $S^3$ by any other manifold of positive constant sectional curvature (e.g. $S^3/\Gamma$, where $\Gamma$ is a finite rotation group which acts freely on $S^3$).

\item (The cylinder) For all $-\infty < t < 0$, let $M_t$ be the standard cylinder $S^2 \times \R$ with metric $g_t = (1-t) g_0 \oplus h$, where $g_0$ is the standard metric on $S^2$ (rescaled to have constant sectional curvature $1/2$) and $h$ is the standard metric on $\R$.  This is an ancient solution (in fact, a gradient shrinking soliton) starting with an infinitely wide round cylinder shrinking down to a cylinder of unit length; if one continued the Ricci flow to positive times, then the cylinder would shrink to a line at time one.  Note that if we made the cylinder compact by replacing $\R$ by $S^1$, then the solution is still ancient but has become collapsed at large scales, because at large scales $r \gg 1$ and at early times $t \ll -r$ the volume of balls of radius $r$ only grows like $r^2$ instead of the expected $r^3$, despite the fact that the curvature never exceeds $O( 1/r^2 )$ on such balls.  Thus we see that the crucial non-collapsing property obtained by Perelman can be used to significantly reduce the possible geometries available for ancient solutions.  It is also worth remarking that in two spatial dimensions, the cylinder $S^1 \times \R$ is also a gradient shrinking soliton and an ancient solution, but is now collapsed at large scales due to the absence of curvature in this case.  (Also, we require $\kappa$-solutions to be non-flat.)

\item (Cigar-type solitons) Let us temporarily work in two spatial dimensions instead of three.  Then, as first observed by Hamilton \cite{hamil-surf}, there is an ancient $t \mapsto M_t$ which at time $t=0$ is the plane $\R^2$ with the metric $dg_0^2 = \frac{dx^2 + dy^2}{1+x^2+y^2}$ (causing the plane to curve somewhat like one end of a cigar), and which at negative times is a isometric to the same plane with metric $(1-t) g_0$.  This solution, known as the \emph{cigar soliton}, is a gradient steady soliton with non-negative curvature, however it is collapsed (it asymptotically resembles the cylinder $S^1 \times \R$ at infinity, and $S^1$ has zero curvature).   Analogues of the cigar soliton - which resemble $S^2 \times \R$ at infinity, and which thus (by the curvature of $S^2$) have enough curvature at infinity to be non-collapsed - are known to exist in three dimensions, see e.g. \cite{bryant}, \cite{cao-soliton}.  These objects are also similar to the \emph{standard solutions} used for the surgery procedure in Section \ref{surgery-sec}.

\item (A non-example) Consider the product of a two-dimensional cigar soliton and a line $\R$.  This is an ancient solution but is collapsed at large scales, because the cigar is collapsed.  It is vital that this example is ruled out from the analysis, because a high-curvature region of Ricci flow which resembles the product of a cigar and a line cannot be removed by surgery.

\end{itemize}

One of the key results of Perelman is that in fact \emph{all} $\kappa$-solutions essentially resemble one of the above examples in a certain quantitative sense, after rescaling to normalise the scalar curvature $R(x)$ to equal $1$.  More on this later.

The analysis of $\kappa$-solutions is rather lengthy and proceeds in several stages; we now consider each stage in turn.

\subsection{Bounded curvature}

We shall begin with an important result on $\kappa$-solutions due to Perelman \cite[\S 11]{per1} (and elaborated in \cite[\S 9.3]{morgan}, \cite[\S 45]{kleiner}, \cite[\S 6.3]{cao}), namely that all $\kappa$-solutions have bounded curvature at time $t=0$ (and hence at all earlier times, by pointwise monotonicity of the scalar curvature for ancient solutions).  This argument is specific to three spatial dimensions (a simpler version of it also works in two dimensions).  

The argument runs by contradiction; if one can find a sequence of points where the curvature goes to infinity, then by looking at rescaled versions of the geometry around these points (or more precisely, slight perturbations of these points where the curvature is essentially a local maximum), and inspecting the rays connecting these points to a fixed origin, one can arrive (after passing to a subsequence) to a limit where one of the scalar curvatures has vanished; this leads to a factorisation of the limit as the Riemannian product of a two-dimensional ancient Ricci flow and a line.  The two-dimensional ancient Ricci flows can be completely classified by a variety of means; one can repeat the argument just stated to show that curvature must be bounded, and in fact furthermore that the manifold is compact.  Then the rounding arguments of Hamilton show that the two-dimensional manifold must in fact have constant curvature, so is basically a sphere.  In short, the limit is basically a cylinder solution\footnote{There is another possibility, that the 2D manifold looks like $RP^2$ instead of $S^2$, but this case can be eliminated by reducing to oriented solutions (this is the approach in \cite{per1}, \cite{kleiner}, \cite{cao}), or by some additional case analysis involving the soul theorem (this is the approach in \cite{morgan}).}.  Undoing the 
scaling, we thus see that the original $\kappa$-solution contains arbitrarily small regions which are approximate rescaled cylinders (the precise term for this is an \emph{$\eps$-neck}, which we shall return to later).  However the soul theorem (see e.g. \cite{petersen}) and the hypothesis of non-negative curvature can be used to eliminate this possibility (see also \cite{chenzhu} for an alternate argument).

\subsection{Curvature decays spatially slower than scaling}\label{curvsec}

Having just obtained a uniform (in $x$) bound on the scalar curvature $R(t,x)$, we now turn to a lower bound in the asymptotic regime $x \to \infty$ (of course this only makes sense for non-compact solutions).  Naive 
dimensional analysis (see Table \ref{dimtable}) suggests that, for a fixed origin $x_0 \in M$, the scalar 
curvature $R(t,x)$ should decay like $1/d(x,x_0)^2$ as $x \to \infty$.  However, the curvature actually decays slower than this.

\begin{proposition}[Infinite asymptotic scalar curvature ratio]\label{isc}\cite[\S 11]{per1}  If $t \mapsto M_t$ is a non-compact $\kappa$-solution, then for each time $t$ we have
$$\limsup_{x \to \infty} R(t,x) d(x,x_0)^2 = +\infty$$ 
for any $x_0 \in M$ (the exact choice of $x_0$ is not relevant).
\end{proposition}

We remark that a slightly weaker version of this proposition was established earlier by Hamilton \cite{hamil-sing}.
 
\begin{proof}(Sketch; for details see \cite[\S 9.6]{morgan}, \cite[\S 40]{kleiner}, \cite[\S 6.4]{cao})  Let us first rule out the case when the limit superior is finite but non-zero at some time $t$.  Then we can find a sequence of radii $r_k$ going to infinity in which $\sup_{d(x,x_k) \sim r_k} R(t,x) \sim r_k^{-2}$;
in particular, as $\kappa$-solutions are non-collapsed, the volume of the annulus $\{d(x,x_k) \sim r_k\}$ is $\sim r_k^3$. We rescale the annuli $\{ d(x,x_k) \sim r_k \}$ by the $r_k$ and take a limit (using Cheeger-Gromov theory) to obtain a Tits cone (the cone over the sphere at infinity) whose maximum scalar curvature is $\sim 1$.  But the Tits cone has some vanishing sectional curvatures.  It is then possible to apply tensorial maximum principles to the equation \eqref{riem} and conclude in fact that the limit is flat, but this contradicts the non-zero scalar curvature.

Now let us rule out the case when the limit superior is zero.  The above argument now shows that the Tits cone is flat, which makes the sphere at infinity round.  This makes the manifold asymptotically flat in a volume sense; using Bishop-Gromov volume comparison theory and the non-negativity of scalar curvature one sees that this is only possible if the manifold was completely flat, a contradiction.
\end{proof}

\subsection{Volume grows spatially slower than scaling}\label{volsec}

Positive curvature tends to contract volume.  Since we have seen that the curvature is larger than what is predicted by scaling, it is perhaps not then surprising that the volume is smaller than is predicted by scaling\footnote{This is consistent with the non-collapsing nature of $\kappa$-solutions, because non-collapsing requires bounded normalised curvature in order to lower bound the volume, and this hypothesis will not be true for extremely large balls with fixed origin.}:

\begin{proposition}[Vanishing asymptotic volume ratio]\label{vac}\cite[\S 11]{per1}
If $t \mapsto M_t$ is a non-compact $\kappa$-solution, then for each time $t$ and origin $x_0 \in M$ the \emph{asymptotic volume ratio} $\lim_{r \to \infty} \Vol_t(B(x_0,r))/r^3$ (which is independent of $x_0$) vanishes.
\end{proposition}

\begin{proof}(Sketch; for details, see \cite[\S 9.6]{morgan}, \cite[\S 40]{kleiner})  This is an induction on dimension argument; to illustrate the key inductive step let us assume that the analogous claim for two dimensions has already been shown and deduce the three-dimensional case.

Fix $t,x_0$. Using the infinite asymptotic scalar curvature ratio and a ``point-picking'' argument of Hamilton \cite{hamil-sing}, one can find $x_k$ and $r_k \to \infty$
with $r_k = o( d(x_0,x_k) )$ such that $\lim_{k \to \infty} R(t,x_k) r_k^2 = \infty$, and such that
$x_k$ essentially maximises the scalar curvature on the ball $B(x_k,r_k)$.  If $\Vol_t( B(x_k,r_k) ) = o(r_k^3)$ then we can use comparison theory to conclude the argument, so suppose that this is not the case.  Rescaling these balls
so that $R(t,x_k) = 1$, we then see that the rescaled balls have curvature at most $1+o(1)$, radius $\tilde r_k \to \infty$ and volume $\sim \tilde r_k^3$, which by comparison theory implies an injectivity radius of $\gtrsim 1$.  One can then use Hamilton's compactness theorem (declaring $x_k$ as the origin of the $B(x_k,r_k)$, of course) to extract a subsequence of these balls which converge to a manifold of scalar curvature at most $1$ (with equality attained at some origin $x_0$), and non-zero asymptotic volume ratio.  One then inspects the rescaled annuli of this manifold (as was done in the proof of Proposition \ref{isc}) to show that the Tits cone splits as the product of a line with a non-flat two-dimensional manifold with non-zero asymptotic volume ratio.   But this can be used to contradict the induction hypothesis\footnote{There is a mild technical issue in ensuring that this two-dimensional manifold is also a $\kappa$-solution, which we will gloss over here.  See \cite[\S 9.6]{morgan} or \cite[\S 40]{kleiner} for details.}.
\end{proof}

We have now seen some ``asymptotic'' connections between high curvature and small volume.  It turns out that one can modify these arguments (again by exploiting Hamilton's compactness theorem) to also obtain some ``local'' connections.  The precise statements are technical (see \cite[\S 6.4]{cao} or \cite[\S 41]{kleiner}), but roughly the gist (allowing for some white lies) is that for any given ball $B(x,r)$ the following statements are morally equivalent:

\begin{itemize}

\item $B(x,r)$ has volume $\gtrsim r^3$.

\item The curvature on any point in $B(x,r)$ is $O(r^{-2})$.

\item The curvature of at least one point in $B(x,r)$ is $O(r^{-2})$.

\end{itemize}

Note that these equivalences already imply Proposition \ref{isc} and Proposition \ref{vac} (and indeed these propositions are used in the proof), since the scalar curvature is strictly positive.  Note that it (morally) doesn't really matter which curvature tensor we use here, since the scalar curvature (essentially) bounds the Riemann and Ricci curvatures thanks to pinching.

There are also a number of useful parabolic estimates which allow one to control derivatives of the Riemann curvature tensor in a pointwise sense by the scalar curvature $R$: indeed we have
$\partial_t^k \nabla_x^l \Riem = O( R^{k + l/2 + 1} )$ for any given $k,l \geq 0$.  

In summary, the scalar curvature $R(x)$ at a point will control the geometry more or less completely out to distances $O(R(x)^{-1/2})$ away from $x$, without any significant collapse of volume or increase in curvature.  Of course $R(x)^{-1/2}$ is the natural length scale associated to $R(x)$ by dimensional analysis.  After that distance, volume collapse will occur, and the curvature will fluctuate, though it has to drop off slower than what scaling predicts.

One consequence of these very strong local controls on curvature (as well as standard analytical tools such as the Rellich compactness theorem) is the following useful compactness theorem of Perelman:

\begin{theorem}[Perelman compactness theorem]\cite[\S 11]{per1} (see also \cite[\S 9.7]{morgan}, \cite[\S 6.4]{cao}, \cite[\S 45]{kleiner}) Let $t \mapsto M_t^k$ be a sequence of $\kappa$-solutions (with uniform control on the non-collapsing), and let $x_k \in M_0^k$ be points. We normalise so that the scalar curvature $R_k( 0, x_k )$ of $(0,x_k)$ on $M_0^k$ is equal to $1$.  Then the sequence $M_t^k$ (with $x_k$ identified as an origin) contains a convergent subsequence whose limit is also a $\kappa$-solution.
\end{theorem}

Very roughly speaking, this compactness theorem gives automatic uniformity on any result involving $\kappa$-solutions, so that once one obtains some control on the geometry of each individual $\kappa$-solution, one automatically gets uniform control on the geometry of all $\kappa$-solutions at once.  We shall therefore gloss over uniformity issues (which we are doing anyway, thanks to the use of
our imprecise $O()$ notation).  

\subsection{Asymptotic gradient shrinking solitons}

It turns out that every $\kappa$-solution has
a (possibly non-unique) gradient shrinking soliton as its ``asymptotic profile'' as $t \to -\infty$.  Roughly speaking, this profile is obtained as follows.  Pick any point $x(0) \in M_0$, and then for any $t < 0$ pick a point $x(t) \in M_t$ which is close to $q$ in reduced length: $l_{0,x(0)}(t,x(t)) \leq 3/2$.  (Such a point can be shown to exist after non-trivial effort by using variational formulae and Hamilton's Harnack inequality to deriving some differential inequalities for parabolic length, and then use standard barrier techniques; see \cite[\S 7.1]{per1}, \cite[\S 23]{kleiner}, \cite[\S 6.2]{cao}, or \cite[\S 7.2]{morgan}.)  In fact there is enough compactness to select $x(t)$ to minimise this reduced length (see \cite[\S 7.2]{morgan}).  The point is that as $t \to -\infty$, the geometry around $(t,x(t))$, shall begin resembling a gradient shrinking soliton after rescaling.

More precisely, let $t_k$ be a sequence of times going to negative infinity.  We rescale the Ricci flow to send $t_k$ to $-1$ and declare $x(t_k)$ to be origin $O$ to obtain a sequence of pointed Ricci flows.  We would like to extract a convergent subsequence from this sequence.  To do this we need uniform control of the geometry of this sequence.  Because we have bounds on the reduced length at the origin $O$, it is possible (thanks to certain differential inequalities for the reduced length, combined with Hamilton's Harnack inequality) to then obtain bounds on the curvature at $O$, and hence by the preceding sections we obtain bounds on curvature at nearby points as well (and also lower bounds on volume, and hence lower bounds on injectivity radius).  This is enough control to apply Hamilton's compactness theorem (or Perelman's compactness theorem) and extract a limiting Ricci flow $t \mapsto M_t^\infty$.  Because each manifold was non-flat, it had reduced volume strictly less than one at time $t=-1$; this reduced volume will be decreasing in $k$ thanks to monotonicity and thus tends to a limit $\tilde V_\infty < 1$.  The reduced lengths $l_k$ for each of the manifolds in the sequence can be shown (with some non-trivial analytic effort) to converge in a moderately strong topology to a limit $l_\infty$ (which morally is the reduced length function of $M^\infty$, though due to various defects in the compactness and initial lack of control on the regularity of the limit flow $M^\infty$, this is not immediately obvious).  We use this pseudo-reduced length $l_\infty$ to define a pseudo-reduced volume $\tilde V_\infty(t)$ on $M^\infty$.  Because the reduced volumes of the rescaled manifolds were converging to a constant $\tilde V_\infty$, one can show by taking limits\footnote{To justify this we need some additional bounds on reduced length away from the origin $O$, but this can be done by using some arguments of Perelman; see \cite{ye}, \cite[\S 6.4]{cao}, or \cite[\S 9.21]{morgan}.} that $\tilde V_\infty(t) = \tilde V_\infty$, thus the pseudo-reduced volume of $\tilde M_t$ is constant in $t$.  Morally speaking, since reduced volume is only constant for gradient shrinking solitons, this should imply that $M_\infty$ is a gradient shrinking soliton.  There are however some annoying analytic issues related to possible
defects of compactness (e.g. the inequality in Fatou's lemma need not be an equality).  Nevertheless it is possible
to extend the above classification of constant reduced volume solutions to constant pseudo-reduced volume by carefully going through the proof of the classification interpreting everything in a distributional sense; see \cite[\S 6.2]{cao}, \cite[\S 38]{kleiner}, or \cite[\S 9.2]{morgan}.  Thus the limit manifold is indeed a gradient shrinking soliton.

\subsection{Classification of gradient shrinking solitons}

Now that we know that $\kappa$-solutions contain rescaled copies of a gradient shrinking soliton as $t \to -\infty$, it is of interest to classify these solitons.  We already know two such examples, the sphere $S^3$ and the cylinder $S^2 \times \R$\footnote{The three-dimensional cigar-type solutions, in contrast, are a gradient \emph{steady} solitons; their asymptotic gradient shrinking soliton is the cylinder $S^2 \times \R$.}.  It turns out that up to some minor details such as quotients by finite group actions, the these are the \emph{only} $\kappa$-solutions which are gradient shrinking solitons.  This important result is due to Perelman \cite[\S 9]{per1} but is mainly proven using earlier technology of Hamilton.  There are several cases to consider in the proof.  If the soliton has at least one zero sectional curvature, it is possible to use Hamilton's maximum principle for the Ricci curvature to conclude that the soliton is a quotient of the cylinder.  If the soliton has positive sectional curvature and is compact, one uses a result of Hamilton \cite{ham1}, which says that all such manifolds become asymptotically round under Ricci flow, to conclude that the gradient-shrinking soliton is perfectly round, and is thus the sphere soliton or a quotient thereof.  The difficult case is when the soliton has strictly positive curvature but is non-compact.  By integrating the equation \eqref{gss} along gradient curves of the potential function $f$ appearing in \eqref{gss}, one can show that the scalar curvature does not go to zero at infinity; standard splitting arguments as in the preceding paragraph then let one show that the manifold is asymptotically a cylinder, with some upper bound on the scalar curvature.  On the other hand, one can show a monotonicity formula for the area of level sets of $f$.  Comparing the lower bounds on these areas against the upper bounds on curvature, one eventually obtains a contradiction to the Gauss-Bonnet theorem, thus excluding this case.  See \cite[\S 6.2]{cao} or \cite[\S 9.4]{morgan} for details.

\subsection{Canonical neighbourhoods}

To summarise so far, we have learnt that the space of $\kappa$-solutions is compact (modulo scaling), and that asymptotically these solutions resemble (after rescaling) either a constant-curvature compact surface (such as a sphere) or a round cylinder.  This strongly suggests that any given $\kappa$-solution at any given time is formed by something resembling the connected sum of rescaled round cylinders and rescaled constant-curvature compact surfaces.  We now (partially) formalise this intuition.

Let us say that a point $x = (t,x)$ in a Ricci flow lives in a \emph{rescaled canonical neighbourhood} if, after rescaling so that $R(x) = 1$, one of the following (somewhat informally phrased) assertions\footnote{For more precise definitions, see \cite[\S 1]{per2}, \cite[\S 9.8]{morgan}, \cite[\S 57]{kleiner}, or \cite[\S 7.1]{cao}.} is true:

\begin{itemize}

\item ($\eps$-neck)  After rescaling, $x$ lies in the center of a region which is within $\eps$ (in a smooth topology such as $C^k$ for suitably large $k$) of the cylinder $[-1/\eps,1/\eps] \times S^2$ with the standard metric, for some suitably small $\eps$.   In some cases, we need to additionally assume that similar statements hold (with slightly wider cylinders) if we evolve backwards in rescaled time for a short (but fixed) duration. This neck models the cylinder solution to Ricci flow.  We also define an \emph{$\eps$-tube} to be a finite or infinite number of $\eps$-necks chained\footnote{One has to be a little careful what ``chaining'' means here; see \cite{morgan} for full details.} together, thus for instance the standard cylinder $\R \times S^2$ is a doubly infinite $\eps$-tube, while cigar-type solitons are a half-infinite $\eps$-tube with a $(C,\eps)$-component (see below) glued on the end.

\item ($C$-component) After rescaling, $x$ lies in a compact manifold diffeomorphic to $S^3$ (or $RP^3$) with diameter $\sim 1$, and all sectional curvatures positive and $\sim 1$.  This component models (rather crudely) the sphere solution to Ricci flow.

\item ($\eps$-round) After rescaling, $x$ lies in a open manifold which is within $\eps$ (in a smooth topology) to a manifold of constant curvature $+1$, such as the sphere $S^3$ or a quotient $S^3/\Gamma$, for some small $\eps$.  This is similar to the $C$-component with more precise control, but possibly with much smaller diameter (because $\Gamma$ could have arbitrarily large).

\item ($(C,\eps)$-cap) After rescaling, $x$ lies in an open region ${\mathcal C}$, diffeomorphic to $S^3 - \{\operatorname{pt}\}$ or $RP^3 - \{\operatorname{pt}\}$ of diameter $O(1)$, scalar curvature positive and $\sim 1$ (with Lipschitz bounds in space and time).  Furthermore, the region ${\mathcal C}$ decomposes into an $\eps$-neck ${\mathcal N}$ (with one end the boundary of ${\mathcal C}$), and a remainder $\overline{Y}$, whose boundary is the other boundary of the $\eps$-neck, and whose interior $Y$ (known as the \emph{core}) contains $x$, and has a non-collapsed geometry at scale $1$.  This cap models cigar-like solutions to Ricci flow.

\end{itemize}

It turns out that every point in a $\kappa$-solution lies in a rescaled canonical neighbourhood.  In fact a more global qualitative description of $\kappa$ solutions is possible:

\begin{theorem}[Description of $\kappa$-solutions]\label{thak} Let $t \mapsto M_t$ be a $\kappa$-solution.  Then one of the following occurs:
\begin{itemize}
\item (Constant-curvature compact) $M_t$ is the round sphere example (or quotient thereof).
\item (Approximately sphere or $RP^3$) $M_t$ is a $C$-component for all $t$.
\item (Cylinder) $M_t$ is the cylinder example, or a quotient thereof by an involution. (In the former case, $M_t$ is diffeomorphic to a punctured $\R^3$.)
\item (Approximate cigar) For each $t$, $M_t$ is diffeomorphic to $\R^3$, has strictly positive sectional curvatures, and is a half-infinite $\eps$-tube with a $(C,\eps)$-cap glued on at one end.
\item (Two-ended cigar) For each $t$, $M_t$ is diffeomorphic to $S^3$ or $RP^3$, has strictly positive sectional curvatures, and is a finite $\eps$-tube with $(C,\eps)$-caps glued on at both ends.
\end{itemize}
\end{theorem}

For more precise formulations of this result, see \cite[\S 11]{per1}, \cite[\S 1]{per2}, \cite[\S 9.8]{morgan}, \cite[\S 58]{kleiner}, \cite[\S 7.1]{cao}.

\begin{proof}(Sketch) There are several cases to consider.  If $M_t$ has vanishing sectional curvature at even just one point and in one direction, then application of Hamilton's maximum principle for Ricci flow allows one to show that $M_t$ (or a two-sheeted covering thereof) factorises as the product of $\R$ with a two-dimensional $\kappa$-solution.  It turns out that there is only one such solution, namely the shrinking two-sphere (this is a consequence of the classification of the sphere $S^2$ as the unique two-dimensional gradient shrinking soliton, combined with Hamilton's observation that any positively curved compact Ricci flow becomes round at the blowup time), and we get the cylinder (or a quotient thereof).

Henceforth we may assume $M_t$ has strictly positive sectional curvatures.  If it is non-compact, then the soul theorem (see e.g. \cite{petersen}) tells us that $M_t$ is diffeomorphic to $\R^3$ and any soul is a point.  Pick one such point, say $x_0$.  Away from this point, one can use the equivalences in Section \ref{volsec} to conclude that the curvature decays slower than scaling, thus $\lim_{x \to \infty} R(t,x) d(x,x_0)^2 = +\infty$ at any time $t$.  Using the classification of asymptotic gradient solitons, one can then conclude that $M_t$ is asymptotically a cylinder at infinity, which turns out to be enough to obtain the approximate cigar case.

One consequence of the above analysis is that every point in a non-compact $\kappa$-solution either lies in a rescaled $\eps$-neck or is in the core of a rescaled $(C,\eps)$-cap; we shall need this result shortly.

The only remaining case is the compact case.  We fix a time $t$ and normalise the maximum scalar curvature to be $+1$.  There are two cases depending on whether the normalised diameter of $M$ is very large ($\gg 1$) for all time or bounded $O(1)$ for at least one time.  In the former case, a rescaling argument exploiting Perelman's compactness theorem and observation in the preceding paragraph shows that every point in $M$ is contained either in the rescaled core of a $(C,\eps)$-cap or a rescaled $\eps$-neck.  Some elementary topological arguments then show that the compact manifold $M$, which is glued together from finitely many of these components, must be homeomorphic to $S^3$, $RP^3$, $RP^3 \# RP^3$, or an $S^2$-fibration over $S^1$.  On the other hand, compact manifolds of positive curvature are known to have finite fundamental group, which rules out the latter two possibilities, and then leads to the two-ended cigar case.

Finally, we have to consider the case when the manifold is compact with bounded normalised diameter at some time $t$.  If at some other time the manifold ceases to have bounded normalised diameter, then by the preceding discussion $M$ is homeomorphic to $S^3$ or $RP^3$ and will then be a rescaled $C$-component (using Perelman's compactness theorem to obtain uniformity in $t$).  So we may assume that we have bounded normalised diameter throughout.
We look at the asymptotic gradient-shrinking soliton.  It must also have bounded normalised diameter, and thus cannot be a cylinder or a quotient; it must be the sphere example.  But then Hamilton's rounding theorem for positively curved compact manifolds implies that $M_t$ is also the round sphere example (or a quotient thereof).  
\end{proof}

\section{Blowup analysis III: approximation by ancient solutions}\label{blow-three}

Having concluded our (lengthy!) analysis of $\kappa$-solutions, we return now to the blowup analysis
of singularities or high-curvature regions of more general solutions to Ricci flow.  Very roughly speaking, Perelman \cite[\S 12]{per1}, \cite[\S 6]{per2} establishes the following three claims regarding a Ricci flow, as one gets sufficiently close to a blow-up time:

\begin{itemize}
\item[(i)] If one rescales a point $(t,x)$ of sufficiently high scalar curvature $R(t,x) \geq K$ so that the scalar curvature is $1$, then a large\footnote{By ``large'' I mean of radius $\gg 1$ in both space and time in rescaled coordinates.  In unrescaled coordinates, this would be a neighbourhood of radius $\gg R(t,x)^{-1/2}$ in space and $\gg R(t,x)^{-1}$ in time.} rescaled neighbourhood of that point is very close to the corresponding neighbourhood of a $\kappa$-solution.
\item[(ii)] If one rescales a point $(t,x)$ of sufficiently high scalar curvature $R(t,x) \geq K$ so that the scalar curvature is $1$, then a large rescaled neighbourhood of that point is a canonical neighbourhood.
\item[(iii)] If one rescales a point $(t,x)$ of sufficiently high scalar curvature $R(t,x) \geq K$ so that the scalar curvature is $1$, then a large rescaled neighbourhood of that point has bounded curvature.
\end{itemize}

For precise statements of these claims see \cite[\S 12]{per1}, \cite[\S 6]{per2}, \cite[\S 10,11]{morgan}, \cite[\S 7]{cao}, \cite[\S 51-52]{kleiner}.  Note that the arguments and statements in \cite[\S 12]{per1} had some inaccuracies, but these were essentially all corrected in \cite[\S 6]{per2}.

The three claims are inter-related.  The claim (i) implies (ii) thanks to Theorem \ref{thak}.
Furthermore, (iii) can be used\footnote{Conversely, one can use the material in Section \ref{volsec} to deduce (iii) from (i), though this direction of implication will not be of use to us.} to imply (i) because of the no-local-collapsing theorem, which parlays the bounded curvature property into an injectivity radius property, which in turn lets\footnote{There is a technical issue here in that one needs to extend the control on the geometry of the rescaled high-curvature region backwards in time somewhat in order to ensure that the limiting solution is indeed ancient.  This subtlety is not fully addressed in Perelman's work, but see \cite[\S 51]{kleiner}, \cite[\S 11.2]{morgan}, \cite[\S 7.1]{cao} for resolutions.} us use Hamilton compactness (and Perelman compactness) to extract the claim.  To conclude the proof of all these claims, Perelman introduces a downward induction on the curvature threshold $K$, essentially reducing matters\footnote{There is the issue of establishing the claims (i)-(iii) for the ``base case'' when $K$ is incredibly large, but this can be dealt with trivially by truncating the Ricci flow just before the first singularity, so that it is still bounded curvature, and then the claims are vacuously true for $K$ large enough.} to showing that (ii) (for a slightly larger value of $K$, say $4K$) implies (iii).  Iteration in $K$ (from incredibly huge down to some moderately large constant) then gives all three claims (i)-(iii).

It remains to show that (ii) for $4K$ implies (iii) for $K$.  Thus, we are assuming that points of really high curvature lie in rescaled canonical neighbourhoods, and want to conclude that points of slightly less high curvature lie in rescaled neighbourhoods without curvature blowup.  Oversimplifying enormously, the argument is as follows.  Suppose for contradiction that the claim failed.  Then one could find a sequence of large neighbourhoods $U_n$ in rescaled Ricci flows where the centre of the neighbourhood had scalar curvature $1$, but that the curvature at other points of the neighbourhood were unbounded.  Furthermore, every point of scalar curvature at least $4$ was contained in a canonical
neighbourhood.  This means that there exists a curve $\gamma_n$ in each $U_n$ which starts at a point of curvature exactly $4$, and passes entirely through canonical neighbourhoods until it reaches a point of arbitrarily large curvature.  It is not hard to show that this scenario cannot happen with the $\eps$-round or $C$-component type of canonical neighbourhoods, so essentially this curve $\gamma$ has to be contained in a long chain of rescaled $\eps$-necks (possibly with one cap on the end, and with different rescaling on each component of the chain, in contrast to an $\eps$-tube; this is a truncated version of an \emph{$\eps$-horn}, which we shall see in the next section).  But the curvature is unbounded, which basically means that the necks are getting narrower and narrower at one end.  It turns out that there is enough control on the geometry inside this chain that one can extract a limit, and obtain an incomplete Ricci flow which is topologically a cylinder $S^2 \times (0,1)$ in space, but has bounded curvature on one end and asymptotically infinite curvature at the other, and consists entirely of $\eps$-necks chained together.  One then approaches the infinite curvature end, rescaling as one goes, and extracts a further limit (using some Aleksandrov-Toponogov theory) to obtain an open non-flat cone (cf. the arguments in Sections \ref{curvsec}, \ref{volsec}).
At this point there is a technical step in which one uses parabolic theory to upgrade the convergence of this limit from
Gromov-Hausdorff convergence to a smoother type of convergence, which we gloss over.  Now we use an argument of Hamilton.  Namely, one observes that 
cones have vanishing sectional curvatures along their rays.  Using Hamilton's maximum principle once again, one can show that the manifold locally splits as a product of a two-dimensional Ricci flow with a line; in particular the scalar curvature stays constant along each ray of the cone.  But this contradicts the fact that the curvature is blowing up as one approaches the apex of the cone.

Details of this argument can be found in \cite[\S 12]{per1}, \cite[\S 4]{per2}, \cite[\S 10-11]{morgan},\cite[\S 51-52]{kleiner}, or \cite[\S 7.1]{cao}.

\section{Horns}\label{hornsec}

We now briefly mention an important technical issue regarding constants.  We have seen that at times close to a singularity, every high-curvature point in a Ricci flow is contained in a canonical neighbourhood, such as an $\eps$-neck.  As it turns out, the $\eps$ parameter here (which is measuring
how close the $\eps$-neck is to a round cylinder) is merely a small absolute constant, such as $10^{-6}$.  For technical reasons later on (having to do with making sure that Ricci-flow-with-surgery has the same type of estimates as smooth Ricci flow) one actually needs to find necks that are much closer to being round, i.e. to find $\delta$-necks for some $\delta \ll \eps$; these very round necks turn out to be excellent locations for performing surgery.  This turns out to be possible by a variant of the arguments used previously.  Very roughly, things proceed as follows.

Let us first pass to the very final time $T_*$ of a smooth Ricci flow, where singularities have already begun to occur.  There will however still be regions of the manifold where the scalar curvature did not pointwise go to infinity; it turns out that in these regions one can still define
a limiting manifold $\Omega$ at time $T_*$, although this manifold will now be incomplete.  On the other hand, we know that the curvature will always go to infinity near any incomplete end of the manifold (in other words, the curvature function is \emph{proper}).  We also know that high curvature points are contained in canonical neighbourhoods.  If these neighbourhoods are $\eps$-round or $C$-components, then they are disconnected from the low-curvature regions of $\Omega$; if the neighbourhood is an $\eps$-cap, and is connected to a low-curvature region then this portion of the manifold is not really close to an incomplete end of the manifold.  Putting all this together, one sees that each incomplete end of $\Omega$ is the infinite end of a half-infinite chain of $\eps$-necks, whose curvature goes to infinity (thus the necks get increasingly narrow at the infinite end).  This shape is known as an \emph{$\eps$-horn}.  Thus the picture of the limiting incomplete manifold $\Omega$ is a large (possibly infinite) number of high curvature connected components (which in addition to $\eps$-round and $C$-component pieces, can also be double $\eps$-horns, a capped $\eps$-horn, or a doubly capped $\eps$-tube), together with the components that contain at least some low-curvature points, which consist of compact sets with finitely many $\eps$-horns attached to them.
See \cite[\S 3]{per2}, \cite[\S 66]{kleiner}, \cite[\S 7.3]{cao}, \cite[\S 11]{morgan} for details.

Now the key point is that $\eps$-horns increasingly resemble a round cylinder as one goes deeper into
the infinite end of the horn.  More precisely, given any $\delta > 0$, it turns out that every point
in the horn of sufficiently high curvature (depending on $\delta$) will be contained in a $\delta$-neck.  The proof of this statement (like many previous arguments here) is by contradiction; if not, then we can find points $x_n$ of increasingly high curvature which are unable to be placed in $\delta$-necks.  Rescaling and taking limits we end up with a solution $M_\infty$ (not necessarily ancient) which has non-negative sectional curvature, is non-collapsed, and is non-flat, and with every point contained in an $\eps$-neck.  By considering what happens to geodesic rays from the $x_n$ to some reference low-curvature point and to a point on the incomplete boundary (i.e. more
Aleksandrov-Toponogov theory) we see that $M_\infty$ has
at least two ends (thus this rules out toroidal examples such as $S^2 \times S^1$).  It turns out that this, combined with negative curvature, already forces $M_\infty$ to be the product of $\R$ with a two-dimensional manifold from Riemannian geometry considerations (see \cite[\S 2.4]{morgan}).  This two-dimensional manifold will of course also obey Ricci flow. Given that $M_\infty$ already has a local structure of an $\eps$-neck, we see that the two-dimensional manifold must be topologically $S^2$ and have positive curvature.  It turns out that one can use
curvature bounds to extend this solution back in time to an ancient solution, at which point the
rounding theory of Hamilton shows that this solution is in fact a round sphere, and so $M_\infty$ is a round cylinder.  In particular the neighbourhoods of $x_n$ must increasingly resemble this cylinder,
a contradiction.  Again, see \cite[\S 3]{per2}, \cite[\S 70]{kleiner}, \cite[\S 7.3]{cao}, \cite[\S 11]{morgan} for details.

\section{Standard solutions}

Now that we have fully analysed the singularities of the Ricci flow, we are almost ready to perform the surgery
procedure that will remove these singularities and replace them with a smooth manifold which can be used to continue the
Ricci flow-with-surgery procedure.  There is however the question of what exactly the smooth manifold will be that will
do this.  This will be the \emph{standard solutions} of Perelman \cite[\S 2]{per2} (the material in which was then
expanded in \cite[\S 7.3]{cao}, \cite[\S 12]{morgan}, or \cite[\S 59-62]{kleiner}).  The advantage of using these
solutions to replace the high curvature regions of a Ricci flow is that their blowup time, and their asymptotic behaviour near that blowup time, is completely understood (to a much greater extent than for arbitrary singularities).

Recall from the previous section that high curvature regions will contain a long $\eps$-tube, which is thus close
to a cylinder.  The idea of surgery is to replace this cylinder-like region with two cigar-type solutions (with
the infinite end of each cigar being one of the two ends of the cylinder); this potentially changes the topology but
will also turn out to take a substantial bite out of the volume, so that this procedure is only performed
finitely often (thanks to the finite extinction result and the bounded growth of volume).

Thus we want a half-infinite solution which resembles a cylinder at one end and is rounded off at the other.  
In light of all the preceding analysis, it is not unreasonable to also request as desirable properties that
the sectional curvatures are all non-negative and that the solution is non-collapsed at all scales $O(1)$.
(It turns out to be neither necessary nor practical to try to make this solution an ancient solution or to enforce 
non-collapsing at large scales.)  This solution will then be referred to as the \emph{standard solution}.  It
is far from unique; its role in the theory is vaguely analogous to the role of a smooth cutoff function in 
harmonic analysis, and so it is the quantitative properties of this solution (such as blowup asymptotics)
which are of importance, rather than the explicit form of the solution.

To construct the standard solution $t \mapsto M_t$ (which will last for time $t \in [0,1)$), we first construct it
at time $t=0$.  What we do here is that we design the initial manifold $M_0$ to be a manifold of revolution
in $\R^4$ (with the Riemannian metric inherited from the ambient Euclidean space $\R^4$), which is 
equal to the half-infinite cylinder $\R^+ \times S^2$ (where $S^2$ is embedded as the unit sphere in $\R^3$ in the
usual manner) with a cap smoothly glued on to the end.  It is not difficult (see \cite[\S 12]{morgan} or \cite[\S 7.3]{cao}) to write down an explicit choice of such a manifold of revolution with non-negative sectional curvature and positive scalar curvature everywhere; this manifold is then isometric with $\R^3$ with some rotationally invariant metric.  We then evolve this manifold by Ricci flow.  Using the local existence and uniqueness theory\footnote{It is possible to avoid use of this theory by approximating the non-compact manifold as a limit of compact manifolds, and exploiting the rotational symmetry assumption heavily to reduce the PDE analysis to what is essentially a $1+1$-dimensional dynamics which can then be analysed by fairly elementary (though lengthy) means.  See \cite[\S 2]{per2}, \cite[\S 12]{morgan}, \cite[\S 59-65]{kleiner} for this approach.  The approach in \cite{cao}, which goes through the non-compact local theory of \cite{chenzhu}, does not actually require rotational symmetry assumptions.} for non-compact manifolds (see \cite{chenzhu}), this solution will exist\footnote{It will also be unique, although this fact is not strictly necessary for the argument, as observed in \cite[\S 59]{kleiner} and attributed to Bernhard Leeb.} up to some maximal time $0 < T_* \leq \infty$ of existence and remain rotationally symmetric throughout.  Since
the initial manifold is asymptotically equivalent to the unit cylinder, which is known to become extinct at time $t=1$, a simple limiting argument shows that the standard solution $t \mapsto M_t$ must asymptotically approach the shrinking
cylinder at all times $0 < t < \min(1,T_*)$ and in particular must develop a singularity by time $t=1$ at the latest;
in other words, $T_* \leq 1$.  Also an application of Hamilton's maximal principle (modified suitably to deal with
the non-compactness of the manifold, see \cite{shi}, \cite[\S 12.2]{morgan}, \cite{chenzhu}) shows that the solution
will always have non-negative sectional curvature and positive scalar curvature (and after $t > 0$ the sectional curvature becomes strictly positive).
It is not hard to demonstrate non-collapsing at scales $O(1)$ at time $t=0$ due to the explicit form of the manifold, and then by Perelman's monotonicity formulae we have non-collapsing at all later times also at the same scales.

Now it is possible to show that blowup cannot occur at any time strictly less than $1$ by the following argument.  Suppose for contradiction that $T$ is strictly less than $1$.  Using the theory of Ricci flows on incomplete manifolds as in the previous section, one can show that as $t \to T$, outside of a compact set $K$, the Ricci flow remains close to
that of the shrinking cylinder solution.  Also, from Proposition \ref{nonneg} we know that the metric is shrinking, and thus metric balls are increasing.  One consequence of this is that the volume of a large ball $B(x,C)$ for some large but fixed $C$ cannot collapse to zero as one approaches the blowup time $T$, because this large ball will contain
a big component outside of $K$ at time $t=0$ and hence at all later times by the increasing nature of the metric balls.  Combining this with Bishop-Gromov comparison theory (cf. Section \ref{nlc}) and the non-negative sectional curvature one can conclude the same statement for small balls: the normalised volume of any ball $B(x,r)/r^3$ for $r = O(1)$ cannot collapse to zero as $t \to T$, even if we also let $r \to 0$ as well.  In other words, the volume of balls cannot grow slower than scaling at scales $O(1)$.

On the other hand, because we are assuming $T$ is the maximal time of existence, then there must be arbitrarily high curvature regions near the blowup time.  By Section \ref{blow-three} these high curvature regions resemble rescaled $\kappa$-solutions.   By Section \ref{volsec} this means that the volume of certain rescaled balls grows slower than scaling.  This rescaled balls get increasingly small as the curvature gets higher, and so we contradict the preceding paragraph.  Thus we must have $T=1$.  (For precise details of this argument, see \cite[\S 2]{per2}, \cite[\S 12.6]{morgan}, \cite[\S 61]{kleiner}, or \cite[\S 7.3]{cao}.)  

Once $T=1$, one can repeat the above types of arguments and eventually obtain a bound
$$ R_{\min}(t) \gtrsim \frac{1}{1-t};$$
the basic reason for this is that if $R_\min$ was ever small at some point $x$ compared to $\frac{1}{1-t}$, then the analysis in Section \ref{blow-three} would imply that the curvature remains bounded near $x$ even up to the blowup time $T=1$.  One can then show that the same holds at all other points, because if there was a transition between bounded limiting curvature and unbounded limiting curvature, the above type of analysis shows that there is no volume collapsing
near the first type of point and volume collapsing near the second type of point, a contradiction.  This leads to a limiting metric at time $t=1$, which must have unbounded curvature since $T=1$ is the maximal time of existence, which implies (by the blowup analysis) that the limiting metric contains arbitrarily small necks.  But this turns out to be incompatible with the positive curvature.  (See \cite[\S 2]{per2}, \cite[\S 61-62]{kleiner}, \cite[\S 12.6]{morgan}, or \cite[\S 7.3]{cao} for details.)

Because of this bound, we know that at times close to $t=1$ the entire manifold consists of high-curvature regions.  One can then use the classification of high curvature regions and basic topology (analogous to the arguments in Section \ref{blow-three}) to show that near the blowup time of a standard solution, every point either lives in a $\eps$-neck or $(C,\eps)$-cap; again, see \cite[\S 2]{per2}, \cite[\S 61-62]{kleiner}, \cite[\S 12.6]{morgan}, or \cite[\S 7.3]{cao} for details.

\section{Surgery}\label{surgery-sec}

We can now informally describe how Ricci flow with surgery is performed.

We first normalise the initial manifold $(M_0,g_0)$ so that the curvature is $O(1)$ and that there is no collapsing at scales $O(1)$.  This is enough to ensure that Ricci flow will remain smooth for a fixed
time interval, say $[0,1]$.  We then split the remaining future time interval $[1,\infty)$ into dyadic pieces $[2^n, 2^{n+1})$ and define Ricci flow with surgery with slightly different surgery parameters on each such dyadic time interval (for technical reasons one needs the parameters on each interval to be much more extreme than those on the previous intervals).  In the case of the Poincar\'e conjecture, when $M_0$ is simply connected, the finite time extinction theory from Section \ref{Extinct-sec} (extended to the case of Ricci flows with surgery) ensure that only finitely many of these intervals need to be considered before the flow becomes completely extinct; however for the purposes of proving the full geometrisation conjecture, finite time extinction is not known in general, and so one must then consider all of these dyadic time intervals.

For the rest of this section we work on a single dyadic time interval $[2^n, 2^{n+1})$, assuming that the flow with surgery has somehow been constructed up to time $2^n$ already.  We shall need two \emph{surgery parameters}, a small \emph{canonical neighbourhood} parameter $r_n$, which has the dimension of length, and a small \emph{surgery control parameter} $\delta_n$, which is dimensionless and measures the accuracy of the canonical neighbourhoods.  

If the flow is non-singular on all of $[2^n, 2^{n+1})$ then we simply perform Ricci flow on the entire 
interval.  Now suppose instead that the flow first becomes singular at some time $T_*$ in the interval $[2^n, 2^{n+1})$.  From the theory in Section \ref{hornsec} we know that the high-curvature regions - the regions where the curvature is $\gg 1/(\delta_n r_n^2)$ - are either disconnected from the low-curvature regions, or else lie in one of a finite number of horns.  Furthermore, there exists a \emph{surgery scale parameter} $h_n$, which is a small parameter with the dimension of length (and much smaller than $\delta_n$ or $r_n$) with the property that every point in the Ricci flow\footnote{Strictly speaking, one needs to extend the analysis of Section \ref{hornsec} from smooth Ricci flows to Ricci flows with surgery.  This turns out to be a remarkably lengthy and non-trivial task, and it is here that it becomes crucial that the $\delta_n, r_n$ are chosen extremely small compared to all previous values of $\delta$ and $r$, but we shall gloss over these very important technical issues in order to focus on the broader picture of the argument.} with scalar curvature at least $h_n^{-2}$ in a horn is contained in the centre of a $\delta_n$-neck.

The high-curvature components which are disconnected from the low-curvature regions are easy to classify: they are either a single compact canonical neighbourhood, or the finite union of non-compact canonical neighbourhoods glued together to form a compact manifold.  Some topology chasing then limits the possibilities for the topology to be $S^3$ or $RP^3$ (or a quotient thereof), or a $S^2$ bundle over $S^1$, or a connected sum of two $RP^3$'s.

At the blowup time $T_*$ we perform the following surgery procedure.  All the high-curvature regions which are disconnected from the low-curvature regions are removed.  For each horn, we locate a point $x_0$ deep inside the horn where the scalar curvature is equal to $h_n^{-2}$, and hence lies inside a $\delta_n$-neck.  In particular, the point $x_0$ is contained inside a two-sphere, which separates the singular part of the horn from the boundary connecting the horn to the low-curvature regions.  We then move from this two-sphere deeper into the horn by a rescaled distance $A$ for some large absolute constant $A$ (and thus move by $Ah_n$ in the non-rescaled metric) - this is still well within the $\delta_n$-neck, and so still gives another two-sphere.  We then cut off and discard everything beyond this two-sphere and replace it by a standard solution (rescaled by $h_n$, and also smoothly interpolated to deal with the fact that the $\delta_n$-neck is not a perfect cylinder).  Now the manifold has become smooth again, and we restart the Ricci flow from this point.  (For more precise descriptions of this surgery procedure, see \cite[\S 4]{per1}, \cite[\S 7.3]{cao}, \cite[\S 67, 71, 72]{kleiner}, \cite[\S 13,14,15.4]{morgan}.)

As we will eventually want to deduce the topology of the initial manifold from the topology at a much later point in the Ricci flow with surgery, it is important to know how to recover the topology of a pre-surgery manifold from the post-surgery manifold.  (The smooth portion of the flow will of course
not affect the topology.) So let us reverse the surgery parameter and
see what happens.  Firstly, reversing the surgery adds a lot of new disconnected high-curvature components to the manifold.  Each point in these components is contained in a canonical neighbourhood.  Some topology-chasing then lets us describe the geometry (and hence topology) of these components as one of the following:

\begin{itemize}

\item A $C$-component, diffeomorphic to $S^3$ or $RP^3$.

\item An $\eps$-round component, diffeomorphic to a finite quotient of $S^3$.

\item An $S^2$ fibration over $S^1$.

\item An $\eps$-tube capped off by two caps (and thus diffeomorphic to $S^3$).

\item An $\eps$-horn capped off at one end.

\item A double $\eps$-horn.

\end{itemize}

The first four cases have no singularities and will just add disconnected components with a well-controlled topology to the pre-surgery manifold.  The other two cases will have singularities which will link up with the ends of the $\eps$-horns attached to the low-curvature components of the manifold.  It is also possible that two $\eps$-horns connected to the low-curvature components of the manifold are attached to each other.  Moving backwards a little bit in time to desingularise the pre-surgery manifold, we see that the effect of the horns is either to connect two components of the post-surgery manifold together, to connect sum a capped cylinder (i.e. a $S^3$) to a component (which does nothing to the topology), or to connect sum a torus $S^2 \times S^1$ to a component.  So if there are only finitely many surgeries before total extinction, each connected component of the manifold at any intermediate time is a connected sum of finitely many manifolds which are either finite quotients of $S^3$, $RP^3$, or $S^2$ fibrations over $S^1$.  More generally, one can show that if a post-surgery manifold obeys Thurston's geometrisation conjecture, then so does the pre-surgery manifold.  
(For details, see \cite[\S 3]{per2}, \cite[\S 72]{kleiner}, \cite[\S 7.3]{cao}, \cite[\S 15]{morgan}.  There are minor variations between these papers due to slightly different definitions of surgery.)  

To complete the proof of Theorem \ref{main} it then suffices (oversimplifying somewhat) to establish the following three claims:

\begin{itemize}

\item[(i)] The Ricci flow with surgery becomes extinct in finite time.

\item[(ii)] The surgery times are discrete (and hence finite, thanks to finite time extinction).

\item[(iii)] The characterisation of high-curvature regions, as used in the above analysis, remains valid even in the presence of surgery.  In particular, the manifold remains non-collapsed at all sufficiently small scales.

\end{itemize}

(For the geometrisation conjecture, (i) can fail\footnote{Indeed, if the fundamental group is \emph{not} a free product of finite groups and infinite cyclic groups - which is the case for instance if any prime component of the manifold has a contractible universal covering, such as a hyperbolic $3$-manifold - then the topological control of surgery discussed earlier shows that the Ricci flow cannot become fully extinct in finite time.  We thank John Morgan for clarifying this point.}, but instead one needs to show that the geometrisation conjecture is asymptotically true in the limit $t \to \infty$.)

The discreteness (ii) is the easiest to establish.  Without loss of generality we may work in a single dyadic time interval $[2^n, 2^{n+1})$.  At each surgery time, either a connected component of the manifold is removed, or an $\eps$-horn is replaced with a standard solution at the scale $h_n$.  In the latter case, the total volume of the manifold is reduced by $\gtrsim h_n^3$, and the number of connected components increased by at most one, for each horn.  In the former case, the number of connected components are decreased by at least one.  On the other hand, because the Ricci curvature is bounded from below, we know that the total volume of the manifold during the smooth part of the Ricci flow can only grow exponentially, and thus stays bounded on a compact time interval such as $[2^n, 2^{n+1}]$.  Putting all this together we see that there can only be finitely many surgeries on any compact interval, i.e. that the surgery times are discrete.  See \cite[\S 4.4]{per2}, \cite[\S 17.2]{morgan}, \cite[\S 72]{kleiner}, \cite[\S 7.4]{cao}.

To establish (iii), the key point is to ensure that some analogue of the monotonicity of the Perelman
reduced volume (see Section \ref{perelcrit}) for the Ricci flow, continues to hold for Ricci flow with surgery.  This turns out to be rather technical but straightforward to verify; the main point is to verify that the curves $\gamma$ which nearly attain the parabolic length \eqref{length-func} do not pass through the regions that are cut away by surgery, and so behave (locally) as if one had a smooth Ricci flow.  In order to achieve this, it becomes essential that the surgery parameters on the dyadic
time interval $[2^n, 2^{n+1}]$ are chosen to be sufficiently small compared to the surgery parameters
on previous intervals.  See \cite[\S 5]{per2}, \cite[\S 76-79]{kleiner}, \cite[\S 7.4]{cao}, \cite[\S 16]{morgan} for details.  One technical point that arises for this analysis is that the quantitative
bounds for the non-collapsing on the time interval $[2^n, 2^{n+1}]$ get substantially worse as $n$ increases, although this ultimately does not cause a major difficulty since it simply means that the
surgery parameters on later intervals have to be chosen even smaller.

Claim (i) is proven in \cite[\S 18]{morgan}.  The strategy is to adapt the finite extinction arguments for Ricci flow of Perelman \cite{per3} as sketched in Section \ref{Extinct-sec} to Ricci flow with surgery.  One can begin with some topological arguments which ensure that after finitely many ``exceptional'' surgeries, the surgeries are all among homotopically trivial $2$-spheres (and so no longer affect $\pi_2(M_t)$).  One then observes that the geometric invariant of the minimal energy $W_2(t)$ of a homotopically non-trivial $2$-sphere in $M_t$ is not only decreasing under Ricci flow (as already noted in Section \ref{Extinct-sec}) but is also
non-increasing under surgery.  As a consequence one readily deduces that all homotopically non-trivial $2$-spheres cease to exist in finite time, thus $\pi_2(M_t)$ is trivial for large $t$.

A similar (but lengthier, and considerably more technical) argument also works for the min-max energy $\tilde W_3(t)$ in Section \ref{Extinct-sec} to ensure $\pi_3(M_t)$ also becomes trivial for 
large $t$.  As mentioned in Section \ref{Extinct-sec}, there is a significant new technical issue unrelated to surgery arising from the need to understand curve shortening flow in order to demonstrate the monotonicity of $\tilde W_3(t)$ even in the absence of surgery.  Once $\pi_3(M_t)$ is trivial one can use topological arguments (combined with the existing control on $\pi_1(M_t)$ and $\pi_2(M_t)$ to force that $M_t$ is empty, and thus we have finite time extinction as desired.  This finally establishes the Poincar\'e conjecture!

\end{document}